\documentclass[11pt,centertags,oneside]{amsart}

\usepackage{amsmath,amstext,amsthm,amssymb,amsfonts,amscd}
\usepackage{mathrsfs}
\usepackage[colorlinks=true]{hyperref}
\usepackage{graphicx}
\usepackage[T1]{fontenc}
\usepackage[latin1]{inputenc}
\usepackage{color,tocvsec2}
\usepackage{typearea}
\usepackage{tikz-cd}
\usepackage{enumerate}
\usepackage{comment}
\usepackage{lscape}
\usepackage[all]{xy}
\usepackage{fouriernc}

\makeatletter
\newcommand{\smallbullet}{} % for safety
\DeclareRobustCommand\smallbullet{%
  \mathord{\mathpalette\smallbullet@{0.8}}%
}
\newcommand{\smallbullet@}[2]{%
  \vcenter{\hbox{\scalebox{#2}{$\m@th#1\bullet$}}}%
}
\makeatother
\usepackage[a4paper,width=14.66cm,top=2.52cm,bottom=2.52cm]{geometry}

\usepackage{multicol}        % used for the two-column index
\makeindex             % used for the subject index
\theoremstyle{definition}
\newtheorem{lem}{Lemma}[section]
\newtheorem{theorem}[lem]{Theorem}
\newtheorem{proposition}[lem]{Proposition}
\newtheorem{cor}[lem]{Corollary}

\newtheorem{definition}[lem]{Definition}

\theoremstyle{remark}
\newtheorem{remark}[lem]{Remark}

\numberwithin{equation}{section}
\numberwithin{figure}{section}

\newcommand{\cE}{\mathcal{E}}

\newcommand{\Z}{\mathbb{Z}}
\newcommand{\R}{\mathbb{R}}
\newcommand{\C}{\mathbb{C}}

\newcommand{\kg}{\mathfrak{g}}
\newcommand{\kh}{\mathfrak{h}}
\newcommand{\kk}{\mathfrak{k}}

\newcommand{\kn}{\mathfrak{n}}
\newcommand{\kp}{\mathfrak{p}}

\newcommand{\kt}{\mathfrak{t}}

\newcommand{\ka}{\mathfrak{a}}

\DeclareMathOperator{\ind}{\mathrm ind}
\DeclareMathOperator{\Tr}{\mathrm Tr}
\DeclareMathOperator{\Ad}{\mathrm Ad}

\newcommand{\bs}{\backslash}

\renewcommand{\(}{\left(}
\renewcommand{\)}{\right)}

\begin{document}

\title{Lefschetz formula for locally symmetric spaces: the real rank one case}

%\curraddr{}
\thanks {}

\makeatletter
\@namedef{subjclassname@2020}{%
  \textup{2020} Mathematics Subject Classification}
\makeatother
\begin{abstract}
Let $G$ be a semi-simple real Lie group of real rank one and $\Gamma$ be a discrete subgroup in $G$ such that $\Gamma \bs G$ has finite volume.  Our main result is an explicit formula for the Lefschetz number of Hecke operators. 
\end{abstract}

\date{\today}
\author{Yanli Song} 
\dedicatory{}

\maketitle
\tableofcontents
%%%%%%%%%%%%%%%%%%%%%%%%%%%%%%%%%%%%%%%%%%%%%%%%%%%%%%%%%%%%%%%%%%%%%%%%%%%%%%%%%%%%%%%%%%%%%%%%%%%%%%%%%%%%%%%%%%%%%%%%%%%%%%%%%%%%%%%
\section{Introduction}

Locally symmetric spaces arise naturally as quotients of symmetric spaces by discrete groups, playing a fundamental role in diverse areas such as arithmetic geometry, number theory, and mathematical physics. A central problem in this setting is the spectral decomposition of $L^2$-spaces associated with these quotients and the study of Hecke operators, which encode deep arithmetic information and provide a powerful tool for investigating automorphic forms and representation theory.

Let $X = G / K$ be a Riemannian symmetric space of noncompact type, where $G$ is a connected, noncompact, real semisimple Lie group and $K$ is a maximal compact subgroup. Consider a Dirac-type operator $D_{\mu} = D_\mu^+ \oplus D_\mu^-$ acting on sections of a bundle associated with a $K$-module $V_{\mu}$, whose highest weight is $\mu$. If $\Gamma$ is a discrete subgroup of $G$, the locally symmetric space $X_{\Gamma} = \Gamma \backslash G / K$ inherits an induced operator $D_{\mu}^{\Gamma} = D_{\mu}^{\Gamma, + }\oplus D_{\mu}^{\Gamma, -}$, which is elliptic. 

When $X_{\Gamma}$ is compact, the Atiyah-Singer index theorem provides an explicit topological formula for the index of $D_\mu^\Gamma$ in terms of the $K$-character $V_\mu$. More generally, we consider a  Atiyah-Singer Lefschetz fixed point theorem for the Hecke operator. Given any element $\alpha$ belongs to the commensurator of $\Gamma$ in $G$, we define a double coset $\Xi = \Gamma \alpha \Gamma$ and a Hecke operator $T_\Xi$ associated to $\Xi$. If  $X_{\Gamma}$ is compact, the Dirac-type operator $D^\Gamma_\mu$ is Fredholm, in the sense that $\dim  \operatorname{Ker} \left(D_\mu^{\Gamma, +}\right)$ are finite dimensional. The Lefschetz number is defined as 
\begin{align}
\label{eq:lefschetz}
L\left(T_\Xi, D_\mu^\Gamma\right) \colon = \Tr \left(T_\Xi \mid \operatorname{Ker} \left(D_\mu^{\Gamma, +}\right) \right) - \Tr \left(T_\Xi \mid \operatorname{Ker}\left(D_\mu^{\Gamma, -}\right) \right) .
\end{align}
 Moscovici \cite{MR765558} derived an explicit expression:
\begin{equation}
\label{compact lefschetz}
L\left(T_\Xi, D_\mu^\Gamma\right) = \sum_\xi \operatorname{vol}\left(\Gamma_\xi \bs G_\xi\right) \cdot  L\left(\xi, D_\mu\right),	
\end{equation}
where the sum runs over representatives of the finitely many central and elliptic $\Gamma$-conjugacy classes in $\Xi$, and each $L(\xi, D_\mu)$ is a Lefschetz number for the $\xi$-action on $X$, given by a topological formula. In a special case when $\Xi = \Gamma$, (\ref{compact lefschetz}) gives the $L^2$-index formula:
\begin{equation}
\label{L2 index}
\ind(D^\Gamma_\mu) = \operatorname{vol}(\Gamma \backslash G) \cdot \ind_G(D_\mu),
\end{equation}
where $\ind_G(D_\mu)$ denotes the ``$G$-index,'' computed topologically in \cite{Connes82}.

However, many of the most interesting locally symmetric spaces are noncompact and exhibit ``cusps''.  The problem of computing the Lefschetz number becomes complicated considerably by the existence of a continuous spectrum.  A natural assumption in this context is that $X_{\Gamma}$ has finite Riemannian volume and strictly negative sectional curvature. According to Theorem 2.1 in \cite{MR658511}, $\dim  \operatorname{Ker} \left(D_\mu^{\Gamma, \pm}\right)$ are still finite dimensional even though the Dirac-type operators $D_\mu^{\Gamma, \pm}$ may fail to be Fredholm (when $\mu$ is singular). Thus,  Lefschetz numbers $L\left(T_\Xi, D_\mu^\Gamma\right)$ in (\ref{eq:lefschetz}) are well-defined.  The main objective of this paper is to extend the formula (\ref{compact lefschetz}) to the noncompact setting, leading to the following:
\begin{equation}
\label{noncompact lefschetz}
L\left(T_\Xi, D_\mu^\Gamma\right) = \sum_\xi L\left(\xi, D_\mu\right) +  \text{``correction term''},
\end{equation}
where the correction term accounts for the presence of cusps and is generally no longer topological.

For de Rham operators, Arthur \cite{MR1001841} provided a solution to (\ref{noncompact lefschetz}) in the adelic setting. Later, Goresky and MacPherson \cite{MR1470341} developed an alternative approach using intersection cohomology techniques. For signature and Dolbeault operators, a Lefschetz formula was obtained by Stern in \cite{MR1258905}. Additional treatments of this problem can be found in \cite{MR989698, MR725778, MR830663}.

In this paper, we focus on general Dirac operators. A pioneering contribution was made by Barbasch and Moscovici \cite{Barbasch83}, where they explicitly computed the correction term in the index formula for Dirac operators, i.e., the case when $\Xi = \Gamma$. Recently, Barbasch-Moscovici's formula was re-examined by Guo-Hochs-Wang \cite{GHW-1, GHW-2} in the framework of K-theory.  We now outline the key ideas. Consider the spinor Laplacians:
\begin{align}
\Delta^{\Gamma, \pm}_\mu = D_\mu^{\Gamma, \pm} D_\mu^{\Gamma, \mp}
\end{align}
on the locally symmetric space $X_{\Gamma}$. The corresponding heat semigroups satisfy:
\begin{align}
e^{-t \Delta_\mu^{\Gamma,  \pm}} = R^{\Gamma}\left(h_{t, \mu}^{ \pm}\right),
\end{align}
where $h_{t, \mu}^{\pm}$ are heat kernels on the symmetric space $X$. Although the operators $R^{\Gamma}\left(h_{t, \mu}^{\pm}\right)$ are generally not of trace class, their discrete part, $R_d^{\Gamma}\left(h_{t, \mu}^{ \pm}\right)$, are. Barbasch and Moscovici computed the index by applying the trace formula developed by Warner \cite{MR535763} and Osborne-Warner \cite{MR614517} to $R_d^{\Gamma} \left(h_t^{+} - h_t^{-}\right)$. 

To avoid certain combinatorial complications, Barbasch and Moscovici imposed a restrictive assumption that $\Gamma$ is neat (i.e., the group generated by the eigenvalues of $\gamma \in \Gamma$ contains no roots of unity). This assumption significantly simplifies the structure of $\Gamma$-conjugacy classes by ensuring that all semisimple components of mixed elements lie in the (finite) center. As a result, the ``correction terms'' vanish in surprisingly many cases.

However, this assumption is highly restrictive and is nearly impossible to satisfy when studying Hecke operators. In this paper, we instead require only that $\Gamma$ be an arithmetic discrete subgroup of $G$ such that $\Gamma \backslash G$ has finite volume. Following \cite{Barbasch83}, we express the Lefschetz number as:
\begin{align}
L(T_\Xi, D_\mu^\Gamma)=\operatorname{Tr} \left(T_\Xi \circ  R_d^{\Gamma} \left(h_t^{+} - h_t^{-}\right) \right), \quad t>0. 
\end{align}
Instead of relying on the trace formula in \cite{MR535763, MR614517}, we employ the \emph{invariant trace formula}, as developed by Arthur \cite{MR625344, MR928262, MR939691} for the adelic setting and Hoffmann \cite{MR1724290} for the real setting. This approach eliminates the need for intricate heat kernel estimates as $t \to \infty$, simplifying the analysis. Finally, we express the traces of $T_\Xi \circ  R_d^{\Gamma} \left(h_t^{+} - h_t^{-}\right)$ in terms of certain invariant tempered distributions on $G$, whose ``Fourier transforms'' (at least the discrete part) are known. 

 Compared to the trace formula in \cite{MR535763, MR614517}, the invariant trace formula for Hecke operators contains additional terms, some of which are less explicit. Nevertheless, using the Selberg's principle, we show that all the implicit terms in the invariant trace formula vanish, ensuring that the correction terms in (\ref{noncompact lefschetz}) are completely explicit. The final formulas (Theorem \ref{main thm I} and \ref{main thm II}) we obtained involve  constant terms in the Laurent expansion of certain Epstein zeta functions, characters of discrete series representations (their restrictions to both compact and non-compact Cartan subgroups), some necessary Lie algebra data, and an additional residue term for singular case.  In particular,  the correction terms are generally nonzero.

As an application, by taking $\Xi = \Gamma$, the index formula in (\ref{noncompact lefschetz}) leads to a multiplicity formula for discrete series representations in $R_d^\Gamma$ with sufficiently regular Harish-Chandra parameters, generalizing the results of Osborne, Warner, and DeGeorge \cite{MR518337, MR671316}. If one chooses other non-trivial $\Xi$, additional refined information about the multiplicities of representations in $R_d^\Gamma$ can be extracted. By selecting a sufficiently rich collection of Hecke operators, it may be possible to compute the multiplicities of discrete series representations that are not sufficiently regular, or even those of non-tempered representations in $R_d^\Gamma$.

This paper is organized as follows:
\begin{itemize}
    \item Section \ref{sec notation}: We introduce the necessary background on Lie algebras, symmetric spaces, and spectral decomposition for the regular representation.
    
    \item Section \ref{sec: Dirac}: We define the Dirac-type operator on (locally) symmetric spaces and study its associated spinor Laplacian.
    
    \item Section \ref{sec hecke}: We introduce Hecke operators and define the Lefschetz number in this context.

    \item Section \ref{sec classification}: We classify the elements in the double coset $\Xi$.
    
    \item Section \ref{sec main results}: We state the main results of this paper.
    
    \item Section \ref{sec application}: We discuss applications of the Lefschetz formula in computing the multiplicities of representations in the discrete spectrum $R^\Gamma_d$.
    
    \item Section \ref{sec invariant}: We explore various invariant tempered distributions on $G$.
    
    \item Section \ref{sec trace}: We present the trace formulas for Hecke operators and provide proofs for the main results.
\end{itemize}

\textbf{Acknowledgements}: We would like to thank Peter Hochs, Haluk \c{S}eng\"un, Xiang Tang, and Hang Wang for their inspiring discussions. We
owe special thanks to Werner Hoffmann for his patience in explaining his work to us. Our research are partially supported by the National Science Foundation grant DMS-1952557. The authors acknowledge support of the Institut Henri Poincar\'e (UAR 839 CNRS-Sorbonne Universit\'e), and
LabEx CARMIN (ANR-10-LABX-59-01).

%%%%%%%%%%%%%%%%%%%%%%%%%%%%%%%%%%%%%%%%%%%%%%%%%%%%%%%%%%%%%%%%%%%%%%%%%%%%%%%%%%%%%%%%%%%%%%%%%%%%%%%%%%%%%%%%
\section{Notation and preliminary results}
\label{sec notation}
\subsection{Symmetric spaces and spinors}
Let $G$ be a noncompact connected semisimple Lie group with finite center. We shall assume that $G$ has real rank one.  Throughout this paper, it will be additionally supposed that $G$ is simple and is embedded in the simply connected complex analytic group corresponding to the complexification of the Lie algebra $\kg$ of $G$.

 Let $\theta \in \operatorname{Aut}(G)$ be the Cartan involution. Define  $K \subset G$  as the fixed point set of  $\theta$  in  $G$ . Then  $K$  is a connected maximal compact subgroup of  $G$ with Lie algebra $\kk$.  The involution  $\theta$  acts as an automorphism on  $\mathfrak{g}$, so that  $\mathfrak{k}$  is the eigenspace of  $\theta$  associated with the eigenvalue  $+1$, while the eigenspace corresponding to the eigenvalue  $-1$  is denoted  $\mathfrak{p}$. This gives the Cartan decomposition:
\begin{align}
\mathfrak{g} = \mathfrak{k} \oplus \mathfrak{p}.
\end{align}
Let  $B: \mathfrak{g} \times \mathfrak{g} \rightarrow \mathbb{R}$  be a symmetric, nondegenerate, bilinear form that is both  $G$ and  $\theta$-invariant. Assume that  $B$  is positive on  $\mathfrak{p}$  and negative on  $\mathfrak{k}$. Thus, the restriction  $B_{\mid \mathfrak{p}}$  induces a Euclidean metric on  $\mathfrak{p}$. Define the non-compact symmetric space:
\begin{align}
X = G / K.
\end{align}
Since  $K$  acts on $\mathfrak{p}$  via the adjoint action, the tangent bundle of  $X$  is given by $T X = G \times_K \mathfrak{p}$.
Equipped with the scalar product induced by  $B_{\mid \mathfrak{p}}$, the space  $X$  becomes a Riemannian manifold.
We assume that
\begin{equation}
\label{equal rank}
\operatorname{Rank} G = \operatorname{Rank} K,
\end{equation}
which is equivalent to  $G$  having a compact Cartan subgroup.

We denote by $S_\kp$ the basic spin representation of $\operatorname{Spin}(\kp)$. Via the adjoint action of $\kk$ on $\kp$, $S_\kp$ becomes a $\kk_\C$-module. Under the assumption in (\ref{equal rank}), we know that $\dim \kp$ is even, and so the $\operatorname{Spin}(\kp)$-module $S_\kp$ decomposes into two irreducible components:
\begin{align}
S_\kp = S_\kp^{+} \oplus S_\kp^{-}.
\end{align}
The map $\Ad \colon K \rightarrow \mathrm{SO}(\mathfrak{p})$ lifts to 
\begin{align}
\widetilde{\Ad}: \widetilde{K} \rightarrow \operatorname{Spin}(\mathfrak{p}),
\end{align} 
for a double cover $\widetilde{K}$ of $K$. Let $\widehat{K}_{\text {Spin }}$ be the set of irreducible representations $V$ of $\widetilde{K}$ such that $S_\kp \otimes V$ descends to a representation of $K$.  Let $R_{\text {spin }}(K)$ be the free abelian group generated by $\widehat{K}_{\text {spin }}$. In particular, if $G/K$ is spin, then $\widehat{K}_{\text {spin }} = \widehat{K}$.

\subsection{Root system and weight lattices}
 Let $\kt$ be the Lie algebra of the compact Cartan subgroup $T$ of $G$, and define the root system decomposition:
\begin{align}
R(\kg, \kt) = R(\kk, \kt) \sqcup R(\kp, \kt).
\end{align}
We denote by $W_\kg$ and  $W_\kk$ the associated Weyl group respectively. Fixing a choice of positive root system, we write:
\begin{align}
R^+(\kg, \kt) = R^+(\kk, \kt) \sqcup R^+(\kp, \kt),
\end{align}
and define the standard half-sums of positive roots:
\begin{align}
&	\rho_\kg = \frac{1}{2} \sum_{\alpha \in R^+(\kg, \kt)} \alpha, & \rho_\kk = \frac{1}{2} \sum_{\alpha \in R^+(\kk, \kt)} \alpha, &&\rho_\kp = \frac{1}{2} \sum_{\alpha \in R^+(\kp, \kt)} \alpha.
\end{align}
The weight lattice in $i \kt^*$ is given by:
\begin{align}
\Lambda_\kt^* = \left\{ \mu \in i \kt^* \Bigg| \frac{2 \langle \mu, \alpha \rangle}{\langle \alpha, \alpha \rangle} \in \mathbb{Z}, \forall \alpha \in  R(\kg, \kt)  \right\}.
\end{align}
We further define:
\begin{align}
\Lambda_\kk^{*}= \left\{ \mu \in \Lambda_\kt^* \Bigg| \langle \alpha, \mu + \rho_\kk \rangle > 0, \forall \alpha \in R^+(\kk, \kt) \right\}, 
\end{align}
and 
\begin{align}
\Lambda_\kk^{*, \text{spin}} =\left\{ \mu \in i \kt^* \Bigg| \frac{2 \langle \mu + \rho_\kp, \alpha \rangle}{\langle \alpha, \alpha \rangle} \in \mathbb{Z}, \forall \alpha \in  R(\kg, \kt), \text{and }  \langle \alpha, \mu + \rho_\kk \rangle > 0, \forall \alpha \in R^+(\kk, \kt)  \right\}. 
\end{align}
By the Borel-Weil theorem \cite{MR1153249}, there exist one-to-one correspondences:
\begin{align}
\widehat{K}\longleftrightarrow \Lambda_\kk^*, \hspace{1cm} \widehat{K}_{\text {spin }}\longleftrightarrow \Lambda_\kk^{*, \text{spin}}. 
\end{align}
Moreover, we define the set of regular weights:
\begin{align}
\Lambda_\kk^{*,\text{reg}} = \left\{ \mu \in \Lambda_\kt^* \Big| \langle \alpha, \mu + \rho_\kk \rangle > 0, \forall \alpha \in R^+(\kg, \kt) \right\}\subset \Lambda_\kk^*.
\end{align}
\subsection{Iwasawa decomposition}
Choose a subspace  $\mathfrak{a} \subset \mathfrak{p}$  such that  $\mathfrak{a}$  is maximal abelian in  $\mathfrak{p}$, and let
\begin{align}
\kg = \mathfrak{k} \oplus \mathfrak{a} \oplus \mathfrak{n}
\end{align}
be the Iwasawa decomposition of $\kg$. In this paper, we assume that the Lie group  $G$  has real rank one, meaning that
\begin{align}
\dim (\mathfrak{a}) = 1.
\end{align}
Denote by  $A$  and  $N$  the subgroups of  $G$  corresponding to  $\mathfrak{a}$  and  $\mathfrak{n}$, respectively.

The set of roots for the pair $\( \kg, \mathfrak{a} \)$ is denoted by $R(\kg, \mathfrak{a})$, and we define $R^+(\kg, \mathfrak{a})$ as the subset of positive roots relative to a fixed ordering. In this setting,
\begin{align}
R^+(\kg, \mathfrak{a}) = \{ \lambda, 2\lambda \}, \quad \text{ and } \quad \rho_\ka = \frac{1}{2}\dim \kg_\lambda \cdot \lambda + \dim \kg_{2\lambda} \cdot \lambda.
\end{align}
Denote the corresponding root spaces as follows:
\begin{align}
\mathfrak{n}_1 = \mathfrak{g}_\lambda, \quad \mathfrak{n}_2 = \mathfrak{g}_{2\lambda}.
\end{align}
We shall also use the notation:
\begin{itemize}
\item $\kn = \kn_1 \oplus \kn_2$;
\item $N_1 = \exp(\kn_1), \quad N_2 = \exp(\kn_2), \quad N = \exp(\kn)$.
\end{itemize}
In particular,  $N_2$  forms the abelian part of the nilpotent group  $N$.

Define  $M$  (resp.  $M^{\prime}$ ) as the centralizer (resp. normalizer) of  $A$  in  $K$. Then, we set
\begin{align}
P = M A N, \quad W_A = M^{\prime} / M \cong \mathbb{Z}_2.
\end{align}
Next, consider the set of cuspidal parabolic subgroups of  G  up to  $\Gamma$ -conjugacy:
\begin{align}
E(G, \Gamma) = \left\{ \Gamma\text{-conjugacy classes of } \Gamma\text{-cuspidal parabolic subgroups of }  G  \right\}.
\end{align}
There are only finitely many elements in  $E(G, \Gamma)$. Moreover, for any two parabolic subgroups  $P_1, P_2 \in E(G, \Gamma)$  with Langlands decompositions
\begin{align}
P_1 = M_1 A_1 N_1, \quad P_2 = M_2 A_2 N_2,
\end{align}
there exists  $k \in K$  such that
\begin{align}
P_1 = k P_2 k^{-1}, \quad M_1 = k M_2 k^{-1}, \quad A_1 = k A_2 k^{-1}, \quad N_1 = k N_2 k^{-1}.
\end{align}

\subsection{Parabolic induction and intertwining operators}
Let $G$ be a semisimple Lie group of real rank one with parabolic subgroup $P = MAN$. The \textit{principal series representations} of $G$ are parametrized by pairs  
\begin{align}
(\sigma, \nu) \in \widehat{M} \times i\mathfrak{a}^*,
\end{align}
where  
\begin{itemize}
    \item $\sigma$ is an \textit{irreducible unitary representation} of $M$, and  
    \item $\nu$ is an element of the \textit{dual space} of the Lie algebra of $A$.  
\end{itemize}
Let $H_\sigma$ be the Hilbert space on which the representation $\sigma$ of $M$ acts. Since $M$ is \textit{compact} by assumption, $H_\sigma$ is \textit{finite-dimensional}. For such a pair $(\sigma, \nu)$, the induced representation $\operatorname{Ind}_P^G(\sigma, \nu)$ is defined as the Hilbert space of all measurable functions $u: G \rightarrow H_\sigma$ satisfying
\begin{align}
\begin{aligned}
u(x \cdot man) &= e^{-(\rho_\ka+\nu)(a)} \sigma(m)^{-1} u(x),  \quad \text{for all } m \in M, a \in A, n \in N, x \in G, \\
\|u\|^2 &= \int_K \|u(k)\|^2 d k < \infty.
\end{aligned}
\end{align}
The principal series representation $\pi_{\sigma, \nu}$ acts on $\operatorname{Ind}_P^G(\sigma, \nu)$ via:
\begin{align}
(\pi_{\sigma, \nu}(x) u)(y) = u(x^{-1} y), \quad u \in \operatorname{Ind}_P^G(\sigma, \nu), \quad  x, y \in G.
\end{align}
Since each $u \in \operatorname{Ind}_P^G(\sigma, \nu)$ is uniquely determined by its restriction to $K$, the space $\operatorname{Ind}_P^G(\sigma, \nu)$ can be identified with the Hilbert space of all measurable functions $v: K \to H_\sigma$ satisfying:
\begin{align}
\begin{aligned}
v(km) &= \sigma(m)^{-1} v(k), \quad \text{for all } m \in M, k \in K, \\
\|v\|^2 &= \int_K \|v(k)\|^2 d k < \infty.
\end{aligned}
\end{align}
Thus, for a fixed $\sigma \in \widehat{M}$, all the representations $\pi_{\sigma, \nu}$ act on the same Hilbert space. Moreover, they have the same restriction to $K$, which we denote by $\operatorname{Ind}_P^G(\sigma)$.

Let $L = MA$ be a Levi subgroup of $G$.  In fact, $L$ is a Levi component of two distinct parabolic subgroups: $P = MAN$ and its opposite $\overline{P} = MA\overline{N}$, with unipotent radicals $N$ and $\overline{N}$. We briefly review the theory of Knapp-Stein intertwining operators \cite[Chapter VII]{KnappBook}. Since $\dim(\mathfrak{a}) = 1$, the Weyl group $W_A = M'/M$ has order two. Let $w \in W_A$ be the nontrivial element. Clearly, $w\nu = -\nu$ for all $\nu \in i\mathfrak{a}^*$. It was pointed out in \cite[P. 13]{MR535763} that  for all elements in $\widehat{M}$ are \emph{ramified}, that is
\begin{align}
w \sigma = \sigma, \quad \sigma \in \widehat{M}. 
\end{align}
The \textit{intertwining operators} associated with $w$ are defined by
\begin{equation}
\label{intertwining}	
I_P( \sigma, \nu ) : H_{\sigma, \nu} \to H_{\sigma,-\nu}, 
\end{equation}
where
\begin{align}
(I_P( \sigma, \nu ) u)(x) = \int_{\overline{N}} u(x \overline{n}) \, d\overline{n}.
\end{align}
These operators satisfy the identity:
\begin{align}
I_P( \sigma, \nu ) \pi_{\sigma, \nu}(x) = \pi_{\sigma, -\nu}(x) I_P( \sigma, \nu ).
\end{align}
The operator $I_P(\sigma, \nu)$ admits a meromorphic continuation to $\nu \in \mathfrak{a}^*_\mathbb{C}$. If $\nu \in i\mathfrak{a}^*$ and $\nu \neq 0$, then $\pi_{\sigma, \nu}$ is irreducible, and $I_P( \sigma, \nu )$ is invertible. The Plancherel density $\mu( \sigma, \nu )$ is defined by
\begin{align}
\mu( \sigma, \nu ) I_P( \sigma, \nu ) I_P( \sigma, -\nu ) = \mathrm{Id}.
\end{align}
It follows from \cite[Theorem 14.16, Corollary 14.30]{KnappBook} that $\pi_{\sigma, 0}$ is irreducible if and only if $\sigma$ is $W_A$-stable and $\mu(\sigma, 0) =0$. Otherwise, if $\mu(\sigma, 0) \neq 0$, the parabolically induced representation decomposes into two \textit{limits of discrete series representations}:
\begin{align}
\pi_{\sigma, 0} \cong \pi_{\sigma}^+ \oplus \pi_{\sigma}^-.
\end{align}
The \textit{normalized intertwining operator} is given by
\begin{equation}
\label{normal intertwining}	
U_P( \sigma, \nu ) := r^{-1}(\sigma, \nu) \cdot  I_P( \sigma, \nu ), 
\end{equation}
where $r(\sigma, \nu)$ is an analytic function satisfying  
\begin{align}
U_P( \sigma, \nu ) \circ U_P( \sigma, -\nu ) = \operatorname{id}.
\end{align}
In particular, $U_P(\sigma, \nu)$ is unitary for all $\nu \in i\mathfrak{a}^*$.

Consider the decomposition of the right regular representation $R_M^{\Gamma}$ of $M$ on $L^2(Z_{\Gamma} \backslash M)$. We have:
\begin{align}
L^2\left(Z_{\Gamma} \backslash M\right) \cong \sum_{\substack{\sigma \in \widehat{M} \\ \sigma \mid_{Z_{\Gamma}}= \mathrm{id}}}^{\oplus} m_\Gamma(\sigma) \cdot H_\sigma, \quad m_\Gamma(\sigma) \in \Z,
\end{align}
and 
\begin{align}
R_M^{\Gamma} \cong \sum_{\substack{\sigma \in \widehat{M} \\ \sigma \mid_{Z_\Gamma}= \mathrm{id}}}^{\oplus} m_\Gamma(\sigma) \cdot \sigma.
\end{align}
The definitions of parabolically induced representations extend naturally to each $P \in E(G, \Gamma)$ with appropriate modifications in notation. 

For any $\sigma \in \widehat{M}$ and $\nu \in i\mathfrak{a}^*$, we define
\begin{equation}
\label{def induced repn}
\operatorname{Ind}_\Gamma^G\left(\sigma, \nu\right) = \sum_{P \in E(G, \Gamma)}m_\Gamma(\sigma) \cdot \operatorname{Ind}_{P}^G\left(\sigma, \nu\right),
\end{equation}
with
\begin{align}
\pi^\Gamma_{\sigma, \nu} = \sum_{P \in E(G, \Gamma)}m_\Gamma(\sigma)  \cdot \pi_{\sigma, \nu}.
\end{align}
By the preceding discussion, the non-trivial element $w \in W_A$ induces unitary intertwining operators:
\begin{align}
U^\Gamma(\sigma, \nu) := \sum_{P \in E(G, \Gamma)} U_P(\sigma, \nu) : \operatorname{Ind}_\Gamma^G\left(\sigma, \nu\right) \rightarrow \operatorname{Ind}_\Gamma^G\left(\sigma,-\nu\right),
\end{align}
which satisfy the functional equation:
\begin{align}
U^\Gamma(\sigma, \nu) \circ U^\Gamma( \sigma, -\nu)  = \operatorname{id}, \quad \text{for all } (\sigma, \nu)  \in \widehat{M} \times i\mathfrak{a}^*.
\end{align}
In what follows, we often identify $\mathfrak{a}^*$ with $\mathbb{R}$ via the isomorphism 
\begin{align}
\nu \in \mathbb{R} \longmapsto \left(\frac{\nu}{\|\lambda\|} \right) \lambda.
\end{align}
Abusing notation, we will use the same letter for both the functional on $\ka$ and the corresponding real number.

\subsection{Spectral decomposition  of regular representation}
Let  $\Gamma$  be  an arithmetic subgroup of $G$  and 
\begin{align}
\operatorname{vol}(\Gamma \backslash G) < \infty.
\end{align} 
Denote by  $R^\Gamma$  the right regular representation of  G  on the Hilbert space  $L^2(\Gamma \backslash G)$. For any  $f \in C_c(G)$, we define:
\begin{align}
R^\Gamma(f) = \int_G f(x) \cdot R^\Gamma(x) \; d_G x.
\end{align}
The right regular representation  $R^{\Gamma}$  of  $G$  on  $L^2(\Gamma \backslash G)$  decomposes into a discrete and a continuous part:
\begin{align}
L^2(\Gamma \backslash G) = L_d^2(\Gamma \backslash G) \oplus L_c^2(\Gamma \backslash G), \quad R^{\Gamma} = R_d^{\Gamma} \oplus R_c^{\Gamma}.
\end{align}
The discrete spectrum  $R_d^\Gamma$  acts on  $L_d^2(\Gamma \backslash G)$, decomposing as a direct sum of irreducible representations, each appearing with finite multiplicity:
\begin{equation}
\label{equ dis dec}
L_d^2(\Gamma \backslash G) \cong \bigoplus_{\pi \in \widehat{\left(\Gamma \bs G\right)}_d} m_{\Gamma}(\pi) \cdot H_\pi, \quad m_{\Gamma}(\pi) \in \mathbb{Z},
\end{equation}
where  $m_\Gamma(\pi)$  denotes the multiplicity of  $\pi$. On the other hand, the continuous spectrum  $R_c^\Gamma$  acts on  $L_c^2(\Gamma \backslash G)$, forming a direct integral of principal series representations, without irreducible subrepresentations. To be more precise, we have the following spectral  decomposition \cite[Section 3]{MR535763}:

\begin{theorem}
\label{contin dec}
We have the orthogonal decomposition:
\begin{align}
L^2(\Gamma \backslash G) \cong L_d^2(\Gamma \backslash G) \oplus L_c^2(\Gamma \backslash G),
\end{align}
with
\begin{align}
L_d^2(\Gamma \backslash G) = \bigoplus_{\pi \in \widehat{\left(\Gamma \bs G\right)}_d} m_{\Gamma}(\pi) \cdot H_\pi, 
\end{align}
\begin{align}
L_c^2(\Gamma \backslash G) = \sum_{\substack{\sigma \in \widehat{M} \\ \sigma |_{Z_{\Gamma}}= \mathrm{id}}}^{\oplus} \left[\int_{i\ka^*}\operatorname{Ind}_\Gamma^G( \sigma, \nu ) \, d\nu \right]^{\mathbb{Z}_2}.
\end{align}
\end{theorem}
Here $Z_\Gamma$ is the center of $\Gamma$ and $\operatorname{Ind}_\Gamma^G( \sigma, \nu )$ are parabolic induced representation of $G$ defined before. 

\subsection{Harish-Chandra $L^p$-Schwartz Spaces}

We now recall the definition of the Harish-Chandra $L^p$-Schwartz spaces $\mathcal{C}^p(G)$ for $p > 0$. This definition involves two spherical functions on $G$, namely $\sigma$ and $\Theta$. 

For $x \in G$, the function $\sigma(x)$ is defined as the \emph{geodesic distance} between the cosets $K$ and $x K$ in the symmetric space $X = G/K$. To introduce the function $\Theta$, let $\pi_0$ be the unitary representation of $G$ induced by the trivial representation of a minimal parabolic subgroup. The trivial representation of $K$ occurs exactly once in the restriction of $\pi_0$ to $K$, and we denote by $\xi$ a unit vector in $\mathscr{H}(\pi_0)$ that is fixed under $\pi_0|_K$. Then the function $\Theta$ is given by the corresponding matrix coefficient:
\begin{align}
\Theta(x) = \left\langle \pi_0(x) \xi, \xi \right\rangle, \quad x \in G.
\end{align}
For $p > 0$, the space $\mathcal{C}^p(G)$ consists of all functions $f \in C^{\infty}(G)$ such that for all $m \geq 0$ and for all $X^\alpha, Y^\beta \in \mathscr{U}(\kg)$,
\begin{align}
\sup_{x \in G} (1+\sigma(x))^m \Theta(x)^{-2 / p} \left| L(X^\alpha) R(Y^\beta) f(x) \right| < \infty,
\end{align} 
where $L$ and $R$ denote the left and right derivatives. In this paper, we sometime denote by $\mathcal{C}(G)$ for $p =2$.

\section{Elliptic Operators on Symmetric Spaces}
\label{sec: Dirac}
We recall the standard procedure for constructing invariant differential operators on a homogeneous space. 
\subsection{Index of Dirac operators}
Let $V_\mu \in \widehat{K}_{\text{spin}}$ with highest weight $\mu \in \Lambda^{*,\text{spin}}_\kk$. Define
\begin{align}
E_\mu^{\pm} = V_\mu \otimes S_\kp^{\pm}, \quad E_\mu = E_\mu^{+} \oplus E_\mu^{-}.
\end{align}
These define homogeneous vector bundles over $X$:
\begin{align}
\mathscr{E}^\pm_\mu := G \times_K E_\mu^{\pm} \to X, \quad \mathscr{E}_\mu := G \times_K E_\mu \to X.
\end{align}
Over the locally symmetric space $X_\Gamma$, we set
\begin{align}
\cE^\pm_{\Gamma,\mu} := \Gamma \bs G \times_K E_\mu^{\pm} \to X_\Gamma, \quad \mathscr{E}_{\Gamma,\mu} := \Gamma \bs G \times_K E_\mu \to X_\Gamma.
\end{align}
As usual, we identify the space of smooth sections $C^{\infty}( \mathscr{E}_\mu)$ with the $K$-invariant part 
\begin{align}
\left(C^{\infty}(G) \otimes E_\mu\right)^K
\end{align} 
Similarly, the space of square-integrable sections $L^2( \mathscr{E}_\mu)$ is identified with
\begin{align}
\left(L^2(G) \otimes E_\mu\right)^K.
\end{align}
Let $\mathscr{U}(\kg_{\C})$ be the universal enveloping algebra of $\kg_{\C}$ and $\mathscr{Z}(\kg_\C)$ the center of $\mathscr{U}(\kg_{\C})$. Let $\left\{X_1, \ldots, X_{\dim \kp}\right\}$ be an orthonormal basis of $\kp$, and define the Dirac operator:
\begin{equation}
\label{def of D}
D_\mu^{\pm}=\sum_{i=1}^{\dim \kp} X_i  \otimes c(X_i) \otimes \operatorname{Id}_{V_\mu} \in \mathscr{U}(\kg_\C) \otimes \operatorname{Hom}(E_\mu^{\pm}, E_\mu^{\mp}),
\end{equation}
where $c(X)$ denotes Clifford multiplication by $X \in \kp$. These operators are $K$-invariant and hence define two $G$-invariant first-order differential operators:
\begin{align}
D_\mu^{\pm} : C_c^{\infty}(\mathscr{E}_\mu^{\pm}) \to C_c^{\infty}(\mathscr{E}_\mu^{\mp}),
\end{align}
called the \emph{Dirac operators} with coefficients in $V_\mu$. They are elliptic, adjoint to each other, and their associated spinor Laplacians, as shown by Parthasarathy \cite{Parthasarathy72}, satisfy the relation:
\begin{align}
D_\mu^{\mp} D_\mu^{\pm}=-R(\Omega) \otimes \operatorname{Id}_{E_\mu^{\pm}} + \left(\|\mu+\rho_c\|^2 - \|\rho\|^2\right) \operatorname{Id}.
\end{align}
Here, $\Omega \in \mathscr{Z}(\kg_\C)$ is the \emph{Casimir element} of $\kg$, and $\|\cdot\|$ refers to the Cartan-Killing norm.

Now let $\pi$ be a unitary representation of $G$ on a Hilbert space $H_\pi$, and let $H^{\infty}_\pi$ be the subspace of all $C^\infty$-vectors in $H_\pi$. With $D_\mu$ given by (\ref{def of D}), we define the operator $\pi(D_\mu)$ on $\left(H_\pi \otimes E_\mu\right)^K$ with (dense) domain $\left(H^{\infty}_\pi \otimes E_\mu\right)^K$ by:
\begin{align}
\pi(D_\mu) = \sum_{i=1}^{\dim \kp} \pi(X_i) \otimes \operatorname{Id}_{V_\mu}.
\end{align}
For the right regular representation $R^{\Gamma} = R^{\Gamma}_d \oplus R^{\Gamma}_c$ on $L^2(\Gamma \bs G)$, we define the associated Dirac operators:
\begin{equation}
\label{equ Dirac ls}
D^\Gamma_{\mu} = R^{\Gamma}(D_\mu), \quad D^\Gamma_{\mu, d} = R_d^{\Gamma}(D_\mu), \quad D^\Gamma_{\mu, c} = R_c^{\Gamma}(D_\mu).
\end{equation}
These define the \emph{Dirac operator} on the locally symmetric space:
\begin{align}
\left[L^2(\Gamma \bs G) \otimes E_\mu  \right]^K \cong L^2(X_\Gamma, \cE_{\Gamma,\mu}).
\end{align}
According to \cite[Theorem 2.1]{MR658511}, we have:
\begin{align}
\operatorname{Ker} D^\Gamma_\mu =\operatorname{Ker}D^\Gamma_{\mu, d}, \quad \text{and} \quad \dim \operatorname{Ker} D^\Gamma_{\mu} < \infty.
\end{align}
In particular, the \emph{index} of $D^\Gamma_{\mu}$ is well-defined:
\begin{align}
\operatorname{index} D^\Gamma_\mu = \dim \operatorname{Ker} D^{\Gamma, +}_{\mu} - \dim \operatorname{Ker}  D^{\Gamma, -}_{\mu}.
\end{align}
\begin{theorem}
\label{multi thm}
The following formula holds:
\begin{align}
\operatorname{index} D_{\mu}^{\Gamma} =  \sum_{\pi \in \widehat{G}, \chi_\pi=\chi_{\mu+\rho_c}} m_{\Gamma}(\pi) \left(\dim \operatorname{Hom}_K\left(H_\pi, E_\mu^{+}\right)-\dim \operatorname{Hom}_K\left(H_\pi, E_\mu^{-}\right)\right),
\end{align}
where $\chi_\pi: \mathscr{Z}(\kg_\C) \to \mathbb{C}$ is the \emph{infinitesimal character} of $\pi \in \widehat{G}$, and $\chi_\lambda$ denotes the character of $\mathscr{Z}(\kg_\C)$ associated with $\lambda \in \kt_\C^*$.
\begin{proof}
\cite[(1.3.8)]{Barbasch83}	
\end{proof}
\end{theorem}

\subsection{Heat semigroups and spinor laplacians}

We now turn our attention to the heat semigroups associated with the spinor Laplacians
\begin{align}
\Delta^\pm_\mu = D_\mu^\pm D_\mu^\mp
\end{align}
acting on $L^2(\mathscr{E}_\mu^{\pm}) = \left(L^2(G) \otimes E_\mu^{\pm} \right)^K$. For each $t > 0$, the bounded operator
\begin{align}
h_{t, \mu}^{\pm} = e^{-t \Delta_\mu^{\pm}} : L^2(\mathscr{E}_\mu^{\pm}) \to L^2(\mathscr{E}_\mu^{\pm})
\end{align}
is a smoothing pseudo-differential operator that commutes with the representation of $G$ on $L^2(\mathscr{E}_\mu^{\pm})$. Hence, it takes the form
\begin{align}
\left(h_{t, \mu}^{\pm} u\right)(x) = \int_G h_{t, \mu}^{\pm} \left(x^{-1} y\right) u(y) \, dy, \quad u \in \left(L^2(G) \otimes E_\mu^{\pm}\right)^K, \quad x \in G,
\end{align}
where $h_{t, \mu}^{\pm} : G \to \operatorname{End} E_\mu^{\pm}$ is a smooth, square-integrable function satisfying the covariance property:
\begin{align}
h_{t, \mu}^{\pm}(x) = k \circ h_{t, \mu}^{\pm} \left(k^{-1} x k'\right) \circ k'^{-1}, \quad x \in G, \quad k, k' \in K.
\end{align}

\begin{theorem}
\label{pi ht thm}
For any irreducible unitary representation $\pi$ of $G$, one has
\begin{align}
\pi\left(h_{t, \mu}^{\pm}\right) = Q_\pi^\pm \circ e^{t\left(\pi(\Omega) - \|\mu+\rho_\kk\|^2 + \|\rho_\kg\|^2\right)} \circ Q_\pi^\pm,
\end{align}
where $Q_\pi^\pm$ denotes the projection of $H_\pi \otimes E_\mu^{\pm}$ onto its finite-dimensional $K$-invariant part $\left(H_\pi \otimes E_\mu^{\pm}\right)^K$.
\begin{proof}
\cite[Corollary 2.2]{Barbasch83}.
\end{proof}
\end{theorem}

The following estimate on the growth of the heat kernel $p_t$ is valid for any connected unimodular Lie group.

\begin{proposition}
Let $t > 0$. Then for all $p > 0$,
\begin{align}
h_{t, \mu}^{\pm} \in \left(\mathscr{C}^p(G) \otimes \operatorname{End} E_\mu^{\pm}\right)^{K \times K}.
\end{align}
\begin{proof}
\cite[Proposition 2.4]{Barbasch83}.
\end{proof}
\end{proposition}

Along with the spinor heat kernels $h_{t, \mu}^{\pm}$, we introduce the one-parameter family of functions
\begin{equation}
\label{equ ht}	
h_{t, \mu}(x) = \Tr^{E_\mu^+} \left( h_{t, \mu}^{+}(x) \right) - \Tr^{E_\mu^-} \left( h_{t, \mu}^{-}(x)\right), \quad x \in G, \quad t > 0.
\end{equation}

\begin{theorem}
\label{thm heat com}
The operators $R_d^{\Gamma}(h_{t, \mu})$ are of trace class, and one has 
\begin{align}
\operatorname{index}D_{\mu}^\Gamma = \operatorname{Tr} R_d^{\Gamma}(h_{t, \mu})
\end{align}
for any $t > 0$.	
\begin{proof}
\cite[Proposition 3.2]{Barbasch83}.
\end{proof}
\end{theorem}

\section{Hecke operators and Lefschetz numbers}
\label{sec hecke}
\subsection{Hecke operators}
Let $\Gamma$ be an arithmetic subgroup of $G$ as before. The commensurator is defined as:
\begin{align}
C_G(\Gamma) := \{ g \in G : \Gamma \cap g\Gamma g^{-1} \text{ has finite index in both } \Gamma \text{ and } g\Gamma g^{-1} \}.
\end{align}
Fix an element $\alpha \in C_G(\Gamma)$. The double coset $\Gamma \alpha \Gamma$ decomposes as a disjoint union of $\Gamma$-cosets. More precisely, there exist elements $\alpha_1, \ldots, \alpha_N \in G$ such that 
\begin{equation}
\label{double co dec}
 \Xi := \Gamma \alpha \Gamma=\bigsqcup_{i=1}^N \alpha_i \Gamma, \quad \text{ (disjoint union)}.	
\end{equation}
Equivalently, we have
\begin{equation}
\label{double co dec II}	
\Gamma=\bigsqcup_{i=1}^N \alpha_i \alpha^{-1} \left(\Gamma \cap \alpha \Gamma \alpha^{-1} \right),
\end{equation}
where $N=[\Gamma : \Gamma \cap \alpha \Gamma \alpha^{-1} ]$.

\begin{definition}
Let $f \in C_c(\Gamma \backslash G)$. Considering $f$ as a $\Gamma$-invariant function on $G$, we define the \emph{Hecke operator} associated with the double coset $\Xi$ by
\begin{align}
T_\Xi(f)( x) =\sum_{i=1}^N f\left(\alpha_i^{-1} x\right), \quad x \in  G,
\end{align}
where the elements $\alpha_i$ are as in (\ref{double co dec}). One can verify that $T_\Xi(f) \in  C_c(\Gamma \backslash G)$, and thus $T_\Xi$ extends to a bounded operator:
\begin{align}
T_\Xi : L^2(\Gamma \bs G) \to  L^2(\Gamma \bs G).
\end{align}
In particular, if $\alpha \in \Gamma$, then $\Xi = \Gamma$ and $T_\Gamma = \operatorname{id}$.
\end{definition}

One can check from the definition that 
\begin{align}
T_\Xi \circ R^\Gamma (x) = R^\Gamma(x)\circ T_\Xi, \quad x \in G.
\end{align}
Thus, $T_\Xi$ defines an \textit{intertwining operator}. We define:
\begin{align}
T_\Xi^{ \pm}=T_\Xi \otimes \operatorname{id} \Big|_{\left(L^2(\Gamma \backslash G) \otimes E^{ \pm}_\mu\right)^K},
\end{align}
which makes sense because $T_\Xi$ commutes with $R_{\Gamma}$. Moreover,
\begin{align}
D^{\Gamma, \pm}_\mu \circ T_\Xi^{ \pm}=T_\Xi^{\mp}\circ D^{\Gamma, \mp}_\mu,
\end{align}
so that $T_\Xi^{ \pm}$ preserves the kernel of $D^{\Gamma, \pm}_\mu$.

\subsection{Lefschetz numbers}

\begin{definition}
The \emph{Lefschetz number} of $T_\Xi$ with respect to $D^{\Gamma, +}_\mu$ is defined as:
\begin{equation}
\label{def lefchetz}
L\left(T_\Xi, D^{\Gamma}_\mu\right) :=  \Tr \left( T^+_\Xi \Big|_{\operatorname{Ker}(D_\mu^{\Gamma, +})} - T^-_\Xi \Big|_{\operatorname{Ker}(D_\mu^{\Gamma,-})}  \right) \in \mathbb{C}.
\end{equation}
In particular, 
\begin{align}
L\left(T_\Gamma, D^{\Gamma}_\mu\right) = \operatorname{index}D_{\mu}^\Gamma.
\end{align}
\end{definition}

We now establish some fundamental properties of these Lefschetz numbers. Given $\pi \in \widehat{G}$, let $L_\pi^2(\Gamma \backslash G)$ denote the corresponding isotypical component of $L^2(\Gamma \backslash G)$, and let $R^{\Gamma}_\pi$ denote the restriction of $R^\Gamma$ to $L_\pi^2(\Gamma  \backslash G)$. Then,
\begin{align}
L_d^2(\Gamma  \backslash G)=\sum_{\pi \in \widehat{\left(\Gamma \bs G\right)}_d}^{\oplus} L_\pi^2(\Gamma \backslash G), \quad 
R^{\Gamma}_d=\sum_{\pi \in\widehat{\left(\Gamma \bs G\right)}_d}^{\oplus} R^{\Gamma}_\pi.
\end{align}
Since $T_\Xi$ preserves each isotypical component, we define
\begin{align}
&T_{\Xi, d}^{ \pm}=T_\Xi^{ \pm} \otimes I \Big|_{\left(L_d^2(\Gamma  \backslash G) \otimes E_\mu^{ \pm}\right)^K}, && T_{\Xi, \pi}^{ \pm}=T_\Xi^{ \pm} \otimes I \Big|_{\left(L_\pi^2(\Gamma  \backslash G) \otimes E^{ \pm}_\mu\right)^K},
\end{align}
so that
\begin{align}
T_{\Xi, d}^{ \pm}=\sum_{\pi \in\widehat{\left(\Gamma \bs G\right)}_d}^{\oplus} T_{\Xi, \pi}^{ \pm}.
\end{align}
Moreover,
\begin{align}
\operatorname{Ker} D^{\Gamma, \pm}_\mu=\sum_{\pi \in\widehat{\left(\Gamma \bs G\right)}_d}^{\oplus} \operatorname{Ker} D_{\mu, \pi}^{\Gamma, \pm},
\end{align}
where $D_{\mu, \pi}^{\Gamma, \pm}$ is the restriction of $D_{\mu}^{\Gamma, \pm}$ to $\left(L_\pi^2(\Gamma  \backslash G) \otimes E^{ \pm}_\mu\right)^K$. 
Since $\dim  \operatorname{Ker} D^{\Gamma, \pm}_\mu < \infty$, it follows that
\begin{align}
\operatorname{Ker} D_{\mu, \pi}^{\Gamma, \pm} = 0, \quad \text{for all but finitely many } \pi \in \widehat{\left(\Gamma \bs G\right)}_d.
\end{align}
Thus, we obtain:
\begin{align}
L\left(T_\Xi, D^{\Gamma}_\mu\right)=\sum_{\pi \in\widehat{\left(\Gamma \bs G\right)}_d}\left(\operatorname{Tr}\left(T_{\Xi, \pi}^{+} \big|_{\operatorname{Ker} D_{\mu, \pi}^{\Gamma, +}}\right) -\operatorname{Tr}\left(T_{\Xi, \pi}^{-} \big|_{\operatorname{Ker} D_{\mu, \pi}^{\Gamma, -}}\right)\right).
\end{align}

\begin{theorem}
\label{heat_kernel_lef}
The Lefschetz number can be computed using the heat kernel of the spinor Laplacian:
\begin{align}
L\left(T_\Xi, D^{\Gamma}_\mu\right)=\Tr\left( T_\Xi \circ  R_d^{\Gamma}\left(h_{t, \mu}\right) \right).
\end{align}
for any $t>0$.
\begin{proof}
\cite[(4.3)]{MR765558}.
\end{proof}
\end{theorem}

\section{Classification of elements in $\Xi$}
\label{sec classification}
In this section, we will classify elements in double coset $\Xi = \Gamma \alpha \Gamma$ for some $\alpha \in C_G(\Gamma)$. See details in \cite[Section 7]{MR1001467} or  \cite[Section 5]{MR614517}. 
\subsection{Semisimple elements}
To begin with, recall that a semi-simple element $\xi \in G$ is classified as:
\begin{itemize}
    \item central if $\xi$ belongs to the center $Z_G$ of $G$.
    \item elliptic if $\xi$ is not central but conjugate to an element of $K$.
    \item hyperbolic if $\xi$ is neither central or elliptic. 
\end{itemize}
We decompose the elements in the double coset $\Xi = \Gamma \alpha \Gamma$.
Define:
\begin{align}
\Xi_S := \{ \xi \in \Xi \mid \forall P \in E(G, \Gamma), \, \{\xi\}_\Gamma \cap P = \emptyset \}.
\end{align}
All elements in $\Xi_S$ are semisimple and the centralizer $\Gamma_\xi$ is a uniform lattice in $G_\xi$. The remaining elements form the parabolic part:
\begin{align}
\Xi_P := \Xi - \Xi_S.
\end{align}
Note that $\Xi_P$ may very well contain semi-simple elements. We further refine:
\begin{align}
\Xi_P = \Xi_P(*) \sqcup \Xi_P(**),
\end{align}
where $ \Xi_P(*)$ consists of semisimple elements and   $\Xi_P(**)$ consists of non-semisimple elements.

On the other hand, we have a regular and singular decomposition of $\Xi_P$ as follow  
\begin{itemize}
    \item Regular elements:
    \begin{align}
    \Xi_P(r) := \{ \xi \in \Xi_P \mid N_{\xi} = \{1\}, \text{ and } \xi \in P \}.
    \end{align}
    One can check that $\Gamma_\xi$ fail to be a lattice in $G_\xi$  for any $\xi \in \Xi_P(r)$.
    \item Singular elements:
    \begin{align}
    \Xi_P(s) := \Xi_P - \Xi_P(r).
    \end{align}
\end{itemize}
By \cite[Lemma 7.1]{MR1001467}, all regular elements in $\Xi_P(r)$ are semisimple, i.e.,
\begin{align}
\Xi_P(r) \subseteq \Xi_P(*).
\end{align}
We decompose the singular part:
\begin{align}
\Xi_P(s) = \Xi_P(s, *) \sqcup \Xi_P(s, **).
\end{align}
where:
\begin{itemize}
    \item $\Xi_P(s, *)$ consists of semisimple singular elements. The semi-simple elements in $\Xi_P(s, *)$ have the property that $\Gamma_\xi$ is a non-uniform lattice in $G_\xi$
    \item $\Xi_P(s, **)$ consists of non-semisimple singular elements.
\end{itemize}
Finally, we have 
\begin{align}
\Xi =
\begin{cases}
\text{Semi-simple part}: & \begin{cases}
\text{Uniform lattice in } G_\xi & \Xi_S, \\
\text{Non-uniform lattice in } G_\xi & \Xi_P(s, *), \\
\text{Not a lattice in } G_\xi & \Xi_P(r).
\end{cases} \\
&\\
\text{Non-semi-simple part}: & \Xi_P(s, ** ).
\end{cases}
\end{align}
Furthermore, all semi-simple elements are decomposed as:
\begin{align}
\Xi_S \sqcup \Xi_P(s, *) \sqcup \Xi_P(r) = \Xi_C \sqcup \Xi_E \sqcup \Xi_H,
\end{align}
where:
\begin{itemize}
    \item $\Xi_C$: central elements.
    \item $\Xi_E$: elliptic elements.
    \item $\Xi_H$: hyperbolic elements.
\end{itemize}
By \cite[p. 382]{MR1001467}, $\Xi_E$ and $\Xi_C$ are finite sets.

 \subsection{Non-semi-simple elements}
 For non-semi-simple elements 
 \begin{align}
 \xi \in \Xi_P(s, **),
 \end{align} 
 we have the Jordan decomposition 
 \begin{align}
 \xi = \eta u.
 \end{align} 
where $\eta$ is semi-simiple, $u$ is unipotent, and they commute. Here the centralizer $N_\eta$ of $\eta$ in $N$ is nontrivial. Since $\xi \in \Gamma \cap P$ for some $P \in E(G, \Gamma)$, we know that $u \in N$ and $L$ can be chosen that $\eta \in L$. This Levi component containing $\eta$ is determined up to conjugation by elements of $N_\eta$.

 For the adjoint action of $A$ on $\mathfrak{n}_\eta$, there is a root space decomposition:
\begin{align}
\mathfrak{n}_\eta = \mathfrak{n}_{\eta,1} + \mathfrak{n}_{\eta,2}, \quad \mathfrak{n}_{\eta,1} \subseteq \mathfrak{n}_1, \quad \mathfrak{n}_{\eta,2} \subseteq \mathfrak{n}_2.
\end{align}
Note that $\mathfrak{n}_2$ is abelian. A \emph{Lie triple} is defined as a triple $(X, H, Y) \neq (0,0,0)$ of elements in $\kg$ satisfying the relations:
\begin{align}
[H, X] = 2 X, \quad [H, Y] = -2 Y, \quad [X, Y] = H.
\end{align}
For $\lambda \in \{1, 2\}$, let $\overset{\circ}{\mathfrak{n}}_{\eta,\lambda}$ be the set of elements $X \in \mathfrak{n}_{\eta,\lambda} = \mathfrak{g}_{\eta, \lambda}$ such that there exists $Y \in \mathfrak{g}_{\eta, -\lambda}$ forming a Lie triple:
\begin{equation} \label{lie triple}
(X, H, Y).
\end{equation}
Defining their exponentiated versions:
\begin{align}\overset{\circ}{N}_{\eta,1}= \exp \left(\overset{\circ}{\mathfrak{n}}_{\eta,1} \right), \quad \overset{\circ}{N}_{\eta,2} = \exp \left(\overset{\circ}{\mathfrak{n}}_{\eta,2} \right)
\end{align}
For such a Lie triple, set 
\begin{align}
Z_0 := X - Y \in \mathfrak{k}.
\end{align}
It has been shown in \cite{MR526310} that  conjugacy classes of unipotent elements are one-to-one correspondent to conjugacy classes of Lie triples, as well as conjugacy classes of $Z_0$. By \cite[Lemma 9.1]{MR1001467}, the set of  $\operatorname{Ad}(L_\eta)$-orbits in $\overset{\circ}{\mathfrak{n}}_{\eta,\lambda}$ are finite.

For $t> 0$, define 
\begin{align}
\xi_t=\eta z_t, \quad z_t =  \operatorname{exp} \left(t Z_0\right).
\end{align} 
We can use the family of semi-simple elements $\xi_t$ to approximate $\xi$. 
\begin{proposition}
 There is a real number $t_0$ such that for any $t \leq t_0$, the centralizer $G_{\xi_t}=G_{\xi_0}$, the centralizer $G_{z_t}=G_{Z_0}$ and $G_{\xi_0} \subseteq G_\eta$ where the element 
 \begin{equation}
 \label{equ Z0}
  \xi_0 =\eta \cdot \exp t_0 Z_0.	
 \end{equation}	
 \begin{proof}
 \cite[Proposition 3.4.1]{MR2626602}	
 \end{proof}
\end{proposition}

\subsection{Epstein zeta function}\label{subsection:epstein}

Another key component of the geometric side of the trace formula is the appearance of certain zeta function constants associated with the discrete group $\Gamma$ (or the double coset $\Xi$ in the case of Hecke operators). Our restriction to the real rank one case leads to a significant simplification compared to the general setting considered in \cite{MR903631}. In the terminology of \cite[Section 9]{MR614517}, all elements in the real rank one case are ``orbitally regular''. As a result, the only zeta function that contributes to the description of the unipotent term in the trace formula is the Epstein zeta function.

 Let us briefly review the construction. For details, we refer to \cite{MR903631, MR614517}. For cuspidal subgroups $P$ with Levi component $L = MA$, put
\begin{itemize}
    \item  $\Xi_L := \text{ projection of } \Xi \cap P \text{ along } N$. 
    \item  $\Xi_M := \text{ projection of } \Xi \cap P \text{ along } AN$. Since $M$ is compact, we know that $\Xi_M$ is a finite set. 
\end{itemize}
 The projection $\Xi_L \to \Xi_M$ is finite to one. In particular, $\Gamma_L = \Gamma_M$.  One should keep in mind that $\Xi_L$ and $\Xi_M$ depend on $P$, although this is not reflected in the notation.

 Let $\eta \in \Xi_L$. If $\xi \in \Xi \cap \eta N$ with Jordan decomposition $\xi=s u$, then $s \in\{\eta\}_N$ and $u \in N_s$, thus $\xi \in s N_s$. Introduce
\begin{align}
\mathscr{S}_\eta=\left\{s \in\{\eta\}_N: s \text { is the semi-simple component of some } \xi \in \Xi \cap \eta N\right\} .
\end{align}
Denote by $n_s$ the element in $N$ with $n_ssn_s^{-1} = \eta$. Put
\begin{align}
\Gamma(s)=n_s\left(\Gamma \cap N_s\right)n_s^{-1}, \quad  \Xi(s)=n_s\left(s^{-1} \Xi \cap N_s\right)n_s^{-1}. 
\end{align}
The above sets are invariant under $\left(\Gamma_M\right)_\eta$. Moreover, $\Gamma(s)$ is a uniform lattice in $N_\eta$, and $\Xi(s)$ is a union of  finitely many $\Gamma(s)$-cosets. 

Now suppose that
\begin{align}
\operatorname{dim}\left(N_{\eta, 2}\right)=1
\end{align}
We define \emph{Epstein zeta function} by 
\begin{align}
\zeta_{\eta}^\pm(s, z)=\operatorname{vol}\left(N_\eta / \Gamma(s)\right) \sum_{\{\xi\}_{\left(\Gamma_M\right)_\eta} \subset \Xi(s) \cap N_{\eta, 2}^{ \pm}} \frac{\operatorname{vol}\left(M_{\xi \eta} /\left(\Gamma_M\right)_{\xi \eta}\right)}{\left\|\log \xi\right\|^{\dim \kn_\eta + z}}.
\end{align}
\cite[Lemma 8]{MR903631} shows that  $\zeta^\pm_{\eta}(s, z)$ is absolutely convergent for $\operatorname{Re} z>0$ and admits a meromorphic continuation to the whole complex plane whose only possible singularities are simple poles at 0 and $-1$.  
\begin{definition}
\label{zeta constant}
For any $\eta \in \Xi_L$, we define 
\begin{align}
\zeta_{\eta}^{\pm}(z)=\sum_{\{s\}_{\Gamma \cap N} \subset \mathscr{S}_\eta}\zeta_{\eta}^{\pm}(s, z).
\end{align}
We denote by $C^\pm_\eta$ the constant term of $\zeta_{\eta}^{\pm}(z)$ at $z = 0$. 	
\end{definition}

\section{Main results: the Lefschetz formulas}
\label{sec main results}
Suppose $G$ is a semisimple Lie group of real rank one and $\Gamma$ is an arithmetic subgroup of $G$.  
Let $\Xi = \Gamma \alpha \Gamma$ be a double coset for some $\alpha \in C_G(\Gamma)$. 
Take $\mu \in \Lambda_\kk^{*, \text{spin}}$. We now state the formula for the Lefschetz number:
\begin{align}
 L\left(T_\Xi, D^{\Gamma}_\mu\right).
\end{align}

\begin{lem}
If $\operatorname{Rank} K \neq \operatorname{Rank} G$, then 
\begin{align}
 L\left(T_\Xi, D_\mu^{\Gamma}\right) = 0.
\end{align}
\begin{proof}
See \cite[Theorem 3.4]{MR765558}.
\end{proof}
\end{lem}

This observation eliminates, in an admittedly trivial way, the case when $\operatorname{Rank} K \neq \operatorname{Rank} G$. 
From now on, we assume:
\begin{align}
\operatorname{Rank} K=\operatorname{Rank} G.
\end{align}
The double coset $\Xi$ decomposes as:
\begin{align}
\Xi = \underbrace{\Xi_C \sqcup \Xi_E \sqcup \Xi_H}_{\text{semi-simple}} \sqcup \underbrace{\Xi_P(s, **)}_{\text{non-semisimple}}.
\end{align}
For brevity, we denote $R^+ = R^+(\kg, \kt)$. For any elliptic element $\xi \in \kg$, set:
\begin{align}
R^+(\xi) = R^+ \cap R(\kg_\xi, \kt).
\end{align}

\begin{theorem} \label{main thm I}
If $\mu$ is regular, that is, 
\begin{align}
\left \langle \mu + \rho_\kk, \alpha \right \rangle > 0, \quad \forall \alpha \in R^+(\kg, \kt),
\end{align}
then the Lefschetz number
\begin{align}
 L\left(T_\Xi, D^{\Gamma}_\mu\right) \in \mathbb{C}
\end{align}
is given by the sum of the following terms:

\begin{enumerate}
    \item {\bf{Central contribution}} (see Theorem \ref{central}):
    \begin{align}
    \delta(\mu, \Xi) \cdot \operatorname{vol}(\Gamma \backslash G) \cdot d(\mu+\rho_\kk),
    \end{align}
    where  
    \begin{align}
    \delta(\mu, \Xi) = \begin{cases}
        \#\left[\Xi_C\right] & \text{if } \zeta_{\mu + \rho_\kk}\vert_{\Xi_C} \equiv 1, \\	
        0 & \text{otherwise}.
    \end{cases}
    \end{align}
    Here, $\zeta_{\mu + \rho_\kk}$ is the central character of the discrete series representation with Harish-Chandra parameter $\mu + \rho_\kk$, and $d(\mu+\rho_\kk)$ is the formal degree given in (\ref{formal degree}).

    \item {\bf{Elliptic contribution}} (see Theorem \ref{elliptic I}):
    \begin{equation} \label{equ elliptic}
    \begin{aligned}
        &(-1)^{\frac{\dim \kp}{2}}  \sum_{ \left\{\xi\right\}_\Gamma \in \Xi_E} \operatorname{vol}\left(\Gamma_\xi \backslash G_\xi\right) \cdot d_\xi^{-1} \\
        &\times \sum_{w \in W_\kk/W_{\kk_\xi}} \operatorname{det}(w) \cdot\left(\frac{ \prod_{\alpha \in R^+(\xi)}\left\langle w\left(\mu+ \rho_\kk\right), \alpha\right\rangle \cdot  e^{w\left(\mu+\rho_\kk\right)}(\xi)}{ e^{\rho_\kg}\left(\xi\right)\prod_{\beta \in R^+ \setminus R^+(\xi)} \left(1  - e^{-\beta}\left(\xi\right)\right)}\right),	
    \end{aligned}	
    \end{equation}
    where $d_\xi$ is a constant depending only on the choice of the Haar measure on $G_\xi$ (for one such choice, $d_\xi$ is given in  \cite{MR0219666}).

    \item {\bf{Parabolic Contribution I}} (see Corollary \ref{unipotent}):
    \begin{equation} \label{equ para I}
    \begin{aligned}
        &\sum_{P \in E(G, \Gamma)} \quad \sum_{\{\eta\}_{\Gamma_M} \subset \Xi_{L}(s)} (-1)^{\frac{\dim \kp}{2}} \cdot \delta(\eta) \cdot  \left(c^+_\eta  \cdot C_\eta^+ + c^-_\eta  \cdot C_\eta^-\right)\\
        &\times \sum_{w \in W_\kk} \left( \left[\overline{\left\langle w( \mu+\rho_\kk), Z_0 \right \rangle}\right]^{\frac{\dim \kn_{\eta,1}}{2}}\cdot \prod_{\beta \in R^+(\xi_0) }\left \langle w( \mu+\rho_\kk), \beta\right\rangle \cdot e^{w( \mu+\rho_\kk)}\left(\eta\right) \right). 
    \end{aligned}	
    \end{equation}
     Here the constants $C_\eta^\pm$ can be found in Definition \ref{zeta constant}, $c^\pm_{\eta}$ are defined in Theorem \ref{thm barbasch} and Remark \ref{remark ceta}, and 
	\begin{align}
	 \delta(\eta) = \begin{cases}
	1& \text{if } \kg_\eta \cong \mathfrak{su}(n, 1).  \text{ In this case  } \dim \kn_{\eta, 2} = 1, \text{and we take}\\	
	 &\text{$X$ to be a representative of the non-trivial unipotent }\\
	 &\text{orbit in $\kn_{\eta, 2}$ such that }\left\{X, H, Y\right\} \text{ and } \left\{-X, H, -Y\right\} \\
	 &\text{are two non-conjugated Lie triples, } Z_0 = X- Y\\
	 &\text{and } \xi_0 \text{ is defined in } (\ref{equ Z0});\\
	 &\\
	 0 & \text{otherwise.}
 	\end{cases}
 	\end{align} 
    \item {\bf{Parabolic Contribution II}} (see Theorem \ref{weighted thm}):
    \begin{equation} \label{equ para II}
    \begin{aligned}
        & \frac{(-1)^{\frac{\dim \kp}{2}+1}}{2} \sum_{P \in E(G, \Gamma)}  \quad \sum_{\{\eta\}_{\Gamma_M} \subset\Xi_{L}} \operatorname{vol}\left(M_\eta /\left(\Gamma_M\right)_\eta\right) \cdot \left|\operatorname{det} \operatorname{Ad}_{\kn}(\eta)\right|^{-1 / 2} \\
        &\times \#\left[ \Xi \cap \eta N \colon \Gamma \cap N \right] \cdot \Omega_{\mu + \rho_\kk}\left(\eta\right),
    \end{aligned}
    \end{equation}
    where $\Omega_{\mu + \rho_\kk}$ is defined (\ref{equ dis cons}). 
\end{enumerate}
\end{theorem}

\begin{remark}
In the Lefschetz formula, the central and elliptic contribution are topological (see Theorem \ref{elliptic II}), and the remaining are not. 	
\end{remark}

\begin{theorem} \label{main thm II}
If $\mu$ is singular, meaning that there exists $\alpha \in R(\kp, \kt)$ such that 
\begin{align}
\left \langle \mu+ \rho_\kk, \alpha \right \rangle = 0,
\end{align}
then $\mu + \rho_\kk$ determines a representation $\sigma \in \widehat{M}$ such that the parabolic induced representation decomposes as:
\begin{align}
\pi_{\sigma, 0} \cong \pi_{\sigma}^+ \oplus  \pi_{\sigma}^-.
\end{align}
The Lefschetz number 
\begin{align}
 L\left(T_\Xi, D^{\Gamma}_\mu\right)
\end{align}
equals the sum of:
\begin{enumerate}
    \item {\bf{Elliptic contribution}}: same as in (\ref{equ elliptic}).
    \item {\bf{Parabolic contribution I}}: same as in (\ref{equ para I}).
    \item {\bf{Residue contribution}}:
    \begin{align}
    -\frac{1}{2} \operatorname{Tr}\left(  T_{\Xi, \pi_{\sigma, 0}} \circ U^{\Gamma, +}\left(\sigma,  0\right) \right).
    \end{align}
\end{enumerate}
\end{theorem}

\begin{remark}
In \cite{Barbasch83}, Barbasch and Moscovici work with the case $\Xi = \Gamma$ and made an addition assumption that  $\Gamma$ is neat (i.e., the group generated by the eigenvalues of any $\xi \in \Gamma$ contains no roots of unity). In their seeting, the trace formula simplifies significantly \cite{MR614517}. Specifically, one can verify that:
\begin{align}
\Gamma_L(s) = \Gamma_L = \Gamma_M = \{ e\}, \quad \Gamma_E = \Gamma_L(r) = \emptyset.
\end{align}
Barbasch and Moscovici applied the trace formula obtained in \cite{MR535763} and proved that if $\mu$ is regular, then the index  is given by the sum of the following contributions:

\begin{itemize}
    \item \text{Central contribution:} 
    \begin{align}
    \operatorname{vol}(\Gamma \backslash G) \cdot d(\mu+\rho_\kk).
    \end{align}
    
    \item \text{Parabolic contribution I:} 
    \begin{align}
    \begin{cases}
    C_2(\Gamma) c_n \varepsilon\left(\Psi_\mu\right) \dim V_\mu,  & \text{if } \kg = \mathfrak{su}(2n, 1), \\
    0, & \text{otherwise}.
    \end{cases}
    \end{align}
    
    \item \text{Parabolic contribution II:} 
    \begin{align}
    \begin{cases}
    \frac{(-1)^{\frac{\dim \kp}{2}+1}}{2} \kappa(\Gamma) \varepsilon_0\left(\mu+\rho_\kk\right) 
    \sum_{w \in W_\kk} \operatorname{det}(w) \operatorname{sign} k\left(w\left(\mu+\rho_\kk\right)\right),  
    & \text{if } \kg = \mathfrak{su}(2, 1) \text{ or } \mathfrak{sl}(2, \mathbb{R}), \\
    0, & \text{otherwise}.
    \end{cases}
    \end{align}
\end{itemize}
One can find the definitions of the above constants in \cite[Theorem 7.1]{Barbasch83}. \end{remark}

\section{Application to multiplicities of discrete series}
\label{sec application}
A fundamental problem is the decomposition of the regular representation of $G$ on $L^2(\Gamma \backslash G)$ into irreducible representations. In particular, we aim to compute the multiplicities $m_{\Gamma}(\pi)$ in the following decomposition:
\begin{align}
L_d^2(\Gamma \backslash G) \cong \bigoplus_{\pi \in \widehat{\left(\Gamma \backslash G\right)}_d} m_{\Gamma}(\pi) \cdot H_\pi, \quad m_{\Gamma}(\pi)\in \mathbb{Z}.
\end{align}
The computation of $m_{\Gamma}(\pi)$ has a long history. If $\Gamma \backslash G$ is compact and $\pi$ belongs to the integrable discrete series representations, Langlands \cite{MR212135} obtained a closed finite formula for $m_{\Gamma}(\pi)$. This result for the compact case was further generalized by Hott and Parthasarathy in \cite{MR348041}.

If $\Gamma \backslash G$ is noncompact, the problem becomes considerably more complicated due to the existence of a continuous spectrum. Moreover, the presence of unipotent elements in $\Gamma$ necessitates the consideration of additional distributions, some of which are invariant. Using Selberg's trace formula, Osborne and Warner \cite{MR518337} derived a multiplicity formula for integral discrete series representations. Later, Barbasch and Moscovici \cite{Barbasch83} extended Osborne-Warner's result to all discrete series representations with sufficiently large Harish-Chandra parameters.

To be more precise, by Theorem \ref{multi thm}, we have
\begin{align}
 L\left(T_\Gamma, D^{\Gamma}_\mu\right)= \sum_{\pi \in \Pi_d(\mu)} m_{\Gamma}(\pi) \cdot m(\pi, \mu),
\end{align}
where
\begin{align}
m(\pi, \mu) = \left(\operatorname{dim} \operatorname{Hom}_K\left(H_\pi, E_\mu^{+}\right)-\operatorname{dim} \operatorname{Hom}_K\left(H_\pi, E_\mu^{-}\right)\right),
\end{align}
and $\Pi_d(\mu)$ denotes the packet consisting of the set of irreducible representations with the same infinitesimal character as $\mu +\rho_\kk$.

In \cite{MR696689}, Baldoni Silva and Barbasch computed the number $m(\pi, \mu) \in \mathbb{Z}$ for all $\pi \in \widehat{G}$ and $\mu \in \Lambda^*_\kk$. In particular, if $\mu$ is sufficiently regular, meaning that 
\begin{align}
\left\langle \mu + \rho_\kk - \rho_\kg, \alpha \right \rangle > 0, \quad \alpha \in R^+(\kg, \kt),
\end{align}
then the only class in $\Pi_d(\mu)$ is the discrete series representation $\omega_{\mu + \rho_\kk}$, and  $m(\pi, \mu)  = 1$. Consequently, the formula for $L\left(T_\Gamma, D^{\Gamma}_\mu\right)$ provides a formula for the multiplicity of the discrete series representation with sufficiently regular Harish-Chandra parameters.

The following corollary generalizes \cite[Corollary 7.2]{Barbasch83} and the formula in \cite{MR518337}.
\begin{cor}
Suppose $G$ is a semisimple Lie group of real rank one and $\Gamma$ is an arithmetic, torsion-free subgroup of $G$ such that $\Gamma \backslash G$ has finite volume. For any sufficiently regular $\mu$, the multiplicity of the discrete series $\omega_{\mu+\rho_\kk}$ in $L_d^2(\Gamma \backslash G)$ is the sum of
\begin{align}
\operatorname{vol}(\Gamma \backslash G) \cdot d(\mu+\rho_\kk),
\end{align}
and 
\begin{align*}
        &\sum_{P \in E(G, \Gamma)} \quad \sum_{\{\eta\}_{\Gamma_M} \subset \Xi_{L}(s)} (-1)^{\frac{\dim \kp}{2}} \cdot \delta(\eta) \cdot  \left(c^+_\eta  \cdot C_\eta^+ + c^-_\eta  \cdot C_\eta^-\right)\\
        &\times \sum_{w \in W_\kk} \left( \left[\overline{\left\langle w( \mu+\rho_\kk), Z_0 \right \rangle}\right]^{\frac{\dim \kn_{\eta,1}}{2}}\cdot \prod_{\beta \in R^+(\xi_0) }\left \langle w( \mu+\rho_\kk), \beta\right\rangle \cdot e^{w( \mu+\rho_\kk)}\left(\eta\right) \right),
    \end{align*}
 and 
\begin{align*}
        & \frac{(-1)^{\frac{\dim \kp}{2}+1}}{2} \sum_{P \in E(G, \Gamma)}  \quad \sum_{\{\eta\}_{\Gamma_M} \subset\Xi_{L}} \operatorname{vol}\left(M_\eta /\left(\Gamma_M\right)_\eta\right) \cdot \left|\operatorname{det} \operatorname{Ad}_{\kn}(\eta)\right|^{-1 / 2} \\
        &\times \#\left[ \Xi \cap \eta N \colon \Gamma \cap N \right] \cdot \Omega_{\mu + \rho_\kk}\left(\eta\right).
   \end{align*}	
\end{cor}

If $\mu$ is not sufficiently regular, then the packet $\Pi_d(\mu)$ may contain more than one irreducible representation, some of which may be non-tempered. For example, let $\Gamma$ be a torsion-free, discrete subgroup of $G=SL(2, \mathbb{R})$. Let 
\begin{align}
D_k^{+}, \quad k=1,  2, \ldots
\end{align} 
be the holomorphic discrete series representations of $SL(2, \mathbb{R})$. Then, we have
\begin{align}
 m_{\Gamma}\left(D_k^{+}\right) =  L\left(T_\Gamma, D^{\Gamma}_k\right), \quad k \geq 2
\end{align}
and 
\begin{equation}
\label{exam sl2}	
 m_{\Gamma}\left(D_1^{+}\right) - m_\Gamma \left(\text{trivial representation}\right) =  L\left(T_\Gamma, D^{\Gamma}_1\right).
\end{equation}
This example illustrates that the pattern for $k=1$ differs from that for $k \geq 2$, since $D_1^+$ and the trivial representation have the same Harish-Chandra parameter.

On the spectral side, the Lefschetz number is given by
\begin{align}
 L\left(T_\Xi, D^{\Gamma}_\mu\right) = \sum_{\pi \in \Pi(\mu)} \Tr\left( T_{\Xi, \pi} \right) \cdot m_\Gamma(\pi).
\end{align}
By selecting a sufficiently rich collection of Hecke operators $T_\Xi$, it may be possible to separate the contributions of $\pi \in \Pi(\mu)$ and precisely determine the multiplicity of a single $m_\Gamma(\pi)$.

%%%%%%%%%%%%%%%%%%%%%%%%%%%%%%%%%%%%%%%%%%%%%%%%%%%%%%%%%%%%%%%%%%%%%%%%%%%%%%%%%%%%%%%%%%%%%%%%%%%%%%%%%%%%%%%%%%%%%%%%%%%%%%%%%%%%%%%%%%%%%%%%%%%%%%%%%%%%%%%%%%%%%%%%%%%%%%%%%%%%%%%%%%%%%%%%%%%%%%%%%%%%%%%%%%%%%%%%%%%%%%%%%%%%%%%%%%%%%%%%%%%%%%%%%%%%%%%%%%%%%%%%%%%%%%%%%%%%%%%%%%%%%%%%%%%%%%%%%%%%

\section{Selberg's principle and invariant tempered distribution}
\label{sec invariant}
In this section, we review the Selberg's principle and  examine several invariant tempered distributions that arise in the trace formula.
\subsection{Selberg's principle}

Suppose that $G$ is a semi-simple Lie group of real rank one. We assume that $G$ has a compact Cartan subgroup $T$ with Lie algebra $\kt$. We denote by $\Lambda^*_\kt$ the weight lattice in $\kt^*$ as before.

A tempered distribution $\mathscr{F}$ on $G$ is said to be \emph{invariant} if it satisfies
\begin{align}
\mathscr{F}\left( f^y\right)  = \mathscr{F}\left(f\right), \quad \text{where } f^y(x) = f(yxy^{-1}), \quad x, y \in G.  
\end{align}

\begin{lem}
If $\mathscr{F}$ is an invariant tempered distribution on $G$, then it
\begin{align}
\mathscr{F}\left(h_{t, \mu}\right)
\end{align}
is independent of $t$.
\begin{proof}
It has been proved in \cite[Lemma 5.4]{GHW} that $\mathscr{F}$ defines a continuous trace on $\mathcal{C}(G)$. Since $\mathcal{C}(G)$ is a dense subalgebra in $C^*_r G$ and is closed under holomorphic calculus \cite{MR1914617}, we obtain the isomorphism
\begin{align}
 K_*\left(\mathcal{C}(G)\right) \cong  K_*\left(C_r^* G\right). 
\end{align} 
This proves the lemma.
\end{proof}
\end{lem}

When the Lie group $G$ has real rank one, Arthur \cite{MR2619638} proved that the Fourier transform induces a topological isomorphism from $\mathscr{C}(G)$ onto $\mathscr{C}(\widehat{G}_{\text{temp}})$. Consequently, the Fourier transform of a tempered distribution can be interpreted as a continuous linear functional on $\mathscr{C}(\widehat{G}_{\text{temp}})$.
If $\mathscr{F}$ is an invariant tempered distribution on $G$, its Fourier expansion has the form:
\begin{equation}
\label{Fourier dis}	
\begin{aligned}
	\mathscr{F}(f) = & \underbrace{\sum_{\lambda \in \Lambda^*_\kt} \widehat{\mathscr{F}}(\lambda) \cdot \Theta_\lambda(f)}_{\text{discrete part}} + 
	\underbrace{\sum_{\sigma \in \widehat{M}} \ \int_{-\infty}^{\infty} \Theta_{\sigma, \nu}(f) \cdot \widehat{\mathscr{F}}(\sigma, \nu) \; d \nu}_{\text{continuous part}}\\
	=&	\mathscr{F}^{(d)}(f) + 	\mathscr{F}^{(c)}(f), \quad f\in \mathcal{C}(G).
\end{aligned}
\end{equation}

The Selberg principle is a vanishing result that plays a key role in our computations within the trace formula.

\begin{lem}
\label{lem selberg}
Let $h_{\mu, t}$ be the function defined in (\ref{equ ht}). For any $\sigma \in \widehat{M}$ and $\nu \in i\ka^*$, we have:
\begin{align}
\Theta_{\sigma, \nu} \left(h_{\mu, t}\right) = 0.
\end{align}
\begin{proof}
By Theorem \ref{pi ht thm}, 
\begin{align}
\Theta_{\sigma, \nu}  \left(h^\pm_{\mu, t}\right) = \dim \left(E_\mu^\pm\right) \cdot e^{t\left(\left\|\sigma\right\|^2 - \left\|\nu\right\|^2-\left\|\mu+\rho_\kk\right\|^2+\|\rho_\kg\|^2\right)}
\end{align} 	
Since $\dim(S^+_\kp) = \dim(S^-_\kp)$, the lemma follows immediately.
\end{proof}
\end{lem}

\begin{theorem}[Selberg's Principle]
\label{cor K selberg}
Let $\mathscr{F}$ be an invariant, tempered distribution on $G$. In terms of the Fourier expansion in (\ref{Fourier dis}), we have:
\begin{align}
\mathscr{F} \left(h_{t, \mu}\right) = (-1)^{\frac{\dim \kp}{2}} \cdot \widehat{\mathscr{F}}(\mu + \rho_\kk).
\end{align}
\begin{proof}
By (\ref{Fourier dis}) and Lemma \ref{lem selberg},
\begin{align*}
\mathscr{F}_*^{(c)}\left(h_{t, \mu} \right)=\sum_{\sigma \in \widehat{M}} \ \int_{-\infty}^{\infty} \Theta_{\sigma, \nu}\left(h_{t, \mu}\right) \cdot \widehat{\mathscr{F}}(\sigma, \nu) \; d \nu = 0.
\end{align*}
Thus, 
\begin{align}
\mathscr{F}_*\left(\ind \left(D_\mu\right) \right) = 	\sum_{\lambda \in \Lambda^*_\kt} \widehat{\mathscr{F}}(\lambda) \cdot \Theta_\lambda(f).
\end{align}
On the other hand, by \cite[(3.5)]{{Barbasch83}}, 
\begin{align}
\Theta_\pi \left( h_{t, \mu}\right)=\begin{cases}
\operatorname{dim} \operatorname{Hom}_K\left(H_\pi, E_\mu^{+}\right)-\operatorname{dim} \operatorname{Hom}_K\left(H_\pi, E_\mu^{-}\right)& \text{if } \chi_\pi =\chi_{\mu + \rho_\kk},\\
&\\
0&\text{otherwise}.
\end{cases}
\end{align}
Thus, we conclude:
\begin{align}
\mathscr{F} \left(h_{t, \mu}\right) =  \widehat{\mathscr{F}}(\mu + \rho_\kk) \cdot \left(\operatorname{dim} \operatorname{Hom}_K\left(\omega_{\mu+\rho_\kk}, E_\mu^{+}\right)-\operatorname{dim} \operatorname{Hom}_K\left(\omega_{\mu+\rho_\kk}, E_\mu^{-}\right) \right),
\end{align}
where $\omega_{\mu+\rho_\kk}$ is the discrete series representation with Harish-Chandra parameter $\mu + \rho_\kk$. Finally, from \cite[Section 4]{MR463358}, we have:
\begin{align}
\operatorname{dim} \operatorname{Hom}_K\left(\omega_{\mu+\rho_\kk}, E_\mu^{+}\right)-\operatorname{dim} \operatorname{Hom}_K\left(\omega_{\mu+\rho_\kk}, E_\mu^{-}\right) = (-1)^{\frac{\dim \kp}{2}}.
\end{align}
This completes the proof. 
\end{proof}
\end{theorem}

%%%%%%%%%%%%%%%%%%%%%%%%%%%%%%%%%%%%%%%%%%%%%%%%%%%%%%%%%%%%%%%%%%%%%%%%%%%%%%%%%%%%%%%%%%%%%%%%%%%%%%%%%%%%%%%%%%%%%%%%

\subsection{Discrete series character}

Let $T$ be a compact Cartan subgroup of $G$ and $H$ be a non-compact Cartan subgroup. Let $\omega_{\mu+\rho_\kk}$ denote the discrete series representation with Harish-Chandra parameter $\mu + \rho_\kk$. For each $f \in C_c^{\infty}(G)$, the operator
\begin{align}
 \omega_{\mu+\rho_\kk}(f) = \int_G f(g) \cdot \omega_{\mu+\rho_\kk}(g) \, d g
\end{align}
is a well-defined bounded linear operator. Moreover, $ \omega_{\mu+\rho_\kk}(f)$ is of trace class, and the map
\begin{align}
f \mapsto \Tr \left(\omega_{\mu+\rho_\kk}(f)\right)
\end{align}
defines a distribution on $G$. We denote this distribution by $\Theta_{\mu + \rho_\kk}$ and refer to it as the character of $\omega_{\mu+\rho_\kk}$. According to Harish-Chandra's fundamental results, we have the following properties:

\begin{itemize}
	\item $\Theta_{\mu + \rho_\kk}$ is an invariant tempered distribution on $G$.
	\item $\Theta_{\mu + \rho_\kk}$ is a locally $L^1$ function on $G$ and is real-analytic on the set of regular semisimple elements of $G$, particularly on $T^{\text{reg}}$ and $H^{\text{reg}}$. 
	\item For any $t \in T^{\text{reg}}$, the character is given by:
	\begin{align}
	\Theta_{\mu + \rho_\kk}(t) =\frac{\sum_{w \in W_\kk}\det(w) \cdot  e^{w\left(\mu+\rho_\kk\right)}(t)}{\Delta_T(t)},	
	\end{align}
	where 
	\begin{align}
	\Delta_T(t) = \operatorname{Det}_{\kg/\kt}^{1/2}\left(\operatorname{Id}-\operatorname{Ad}_t\right).
	\end{align}
	\item On $H^{\text{reg}}$, the character $\Theta_{\mu + \rho_\kk}$ is determined as follows. Define:
	\begin{align*}
	A^+ &= \left\{a \in A \mid e^\lambda(a) > 1, \quad \lambda \text{ real positive root}\right\},\\
	A^- &= \left\{a \in A \mid e^\lambda(a) < 1, \quad \lambda \text{ real positive root}\right\}.		
	\end{align*}
	Split $H = H_KA$ with $H_K = H \cap K$, and define:
	\begin{align}
	H^+ = \left(H_K \cdot A^+ \right)\cap H^{\text{reg}}, \quad H^- = \left(H_K \cdot A^- \right)\cap H^{\text{reg}}. 
	\end{align}
	For any $\mu \in \Lambda_\kt^*$, define the sign function:
	\begin{align}
	c(\mu, H^+) = \begin{cases}
 	-1 & \text{if } \mu\left(\sqrt{-1}\left(E_\lambda - E_{-\lambda}\right)\right) > 0,\\
 	+1 & \text{if } \mu\left(\sqrt{-1}\left(E_\lambda - E_{-\lambda}\right)\right) < 0,\\
 	0 & \text{if } \mu\left(\sqrt{-1}\left(E_\lambda - E_{-\lambda}\right)\right) = 0,
 	\end{cases}
	\end{align}
	and define $c(\mu, H^-)  = -c(\mu, H^+)$,
	where $E_{\pm \lambda}$ are the normalized vectors in $\kg_{\pm \lambda}$. Let $\mathbf{c} \colon H \to T$ denote the Cayley transform. Then, for any $h \in H^{\text{reg}}$, we have:
	\begin{equation}
	\Theta_{\mu + \rho_\kk}(h) 	=  \frac{\sum_{w \in W_\kg} \operatorname{det}(w) \cdot c\left(w \mu: J^{ \pm}\right) \cdot   e^{w \mu}(\mathbf{c}(h))}{\Delta_H(h)},
	\end{equation}
	where the sign $\pm$ is chosen based on whether $h \in H^+$ or $H^-$, and 
	\begin{align}
	\Delta_H(h) = \operatorname{Det}_{\kg/\kh}^{1/2}\left(\operatorname{Id}-\operatorname{Ad}_h\right).
	\end{align}
\end{itemize}

\begin{remark}
From the above result, we see that the function
\begin{equation}
\label{equ dis cons}
\begin{aligned}
\Omega_{\mu + \rho_\kk}(h) &:= \Delta_H(h)  \cdot \Theta_{\mu + \rho_\kk}(h)\\
&=\sum_{w \in W_\kg} \operatorname{det}(w) \cdot c\left(w \mu: J^{ \pm}\right) \cdot   e^{w \mu}(\mathbf{c}(h))
\end{aligned}
\end{equation}
extends to a smooth function on $H$.  
\end{remark}

%%%%%%%%%%%%%%%%%%%%%%%%%%%%%%%%%%%%%%%%%%%%%%%%%%%%%%%%%%%%%%%%%%%%%%%%%%%%%%%%%%%%%%%%%%%
\subsection{Semi-simple orbital integral}

Let $\xi \in G$, and let $Z_\xi \subseteq G$ be its centralizer. Let $d(g Z_\xi)$ be the left-invariant measure on $G / Z_\xi $ determined by a Haar measure $d_G$ on $G$. 

\begin{definition}
The \emph{orbital integral} of a measurable function $f$ on $G$ with respect to $\xi$ is given by
\begin{align}
\tau_\xi(f):=\int_{G / Z_\xi} f\left(g \xi g^{-1}\right) d\left(g Z_\xi\right),
\end{align}
provided the integral converges.
\end{definition}

Harish-Chandra proved that the integral converges for $f$ in the Harish-Chandra Schwartz algebra $\mathcal{C}(G)$ \cite[Theorem 6]{MR0219666}. Furthermore, results from \cite{MR320232} show that $\tau_\xi$ defines an invariant tempered distribution on $G$ and is also a measure on $G$. 

A central problem in representation theory is obtaining Fourier inversion formulas in the sense of Harish-Chandra for the distributions $\tau_\xi$, that is, explicit formulas in (\ref{Fourier dis}). This means expressing $\tau_\xi$ as a series of tempered invariant eigen-distributions of $G$. Such formulas have been obtained in \cite{MR0450461} for $\xi$ semisimple and $G$ of real rank one. For $\xi$ regular and $G$ of rank greater than one, formulas are obtained in \cite{MR525674}.

While orbital integrals appear on the geometric side of the trace formula, their spectral side counterparts are characters. The character of a discrete series representation can be realized in the following.

\begin{theorem}[Character Formula for Elliptic Orbital Integral]
\label{elliptic I}
If $\xi$ is elliptic, then 
\begin{align*}
\tau_\xi\left(h_{t, \mu}\right)&=  
 (-1)^{\frac{\dim \kp}{2}}  \cdot d_\xi^{-1} \sum_{w \in W_\kk/W_{\kk_\xi}} \operatorname{det}(w) \\
 &\quad \times \left(\frac{ \prod_{\alpha \in R^+(\xi)}\left\langle w\left(\mu+ \rho_\kk\right), \alpha\right\rangle \cdot  e^{w\left(\mu+\rho_\kk\right)}(\xi)}{ e^{\rho_\kg}(\xi)\prod_{\beta \in R^+ \setminus R^+(\xi)} \left(1  - e^{-\beta}(\xi)\right)}\right),  
\end{align*}
where $d_\xi$ is a constant depending only on the choice of the Haar measure on $G_\xi$. In particular, if $\xi \in T^{\text{reg}}$, then 
\begin{equation}
\label{character orbital}
\tau_\xi\left(h_{t, \mu}\right)= (-1)^{\frac{\dim \kp}{2}} \cdot \Theta_{\mu + \rho_\kk}.
\end{equation}

\begin{proof} 
We provide a sketch; full details can be found in \cite[Theorem 4.14]{MR765558}. By $G$-invariance, we may conjugate $\xi$ so that $\xi \in T$. 

If $\xi \in T^{\text{reg}}$, the formula was proved in \cite[Theorem 3.1]{MR3990784}. For singular $\xi$, the character $\Theta_{\mu + \rho_\kk}$ diverges, but the regularized orbital integral
\begin{align}
\Delta_T(\xi) \cdot \tau_\xi\left(\ind_G\left(D_\mu\right) \right) 
\end{align}
remains continuous on $T$. 

Define the differential operator on $\kt$:
\begin{align}
\Psi_\xi = \prod_{\alpha \in R^+(\xi)} H_\alpha \in \mathscr{U}\left( \kt_\C\right),
\end{align}
where $H_\alpha$ is determined by
\begin{align}
B\left(H, H_\alpha\right) = \alpha(H), \quad \forall H \in \kt.
\end{align}
Harish-Chandra showed in \cite{MR0219666} that 
\begin{align}
\tau_\xi\left(h_{t, \mu}\right) = \frac{e^{\rho_\kg}(\xi)}{ d_\xi \cdot \prod_{\beta \in R^+ \setminus R^+(\xi)}\left(e^\beta(\xi)-1\right)} \cdot \lim_{\substack{\eta \to \xi\\ \eta \in T^{\text{reg}}}} \Psi_\xi \left(\Delta_T(\eta) \cdot \tau_\eta\left(\ind_G\left(D_\mu\right) \right) \right). 
\end{align}
\end{proof}
\end{theorem}

\begin{theorem}[Central Orbital Integral]
\label{central}
If $\xi$ is central, then 
\begin{align}
\tau_\xi\left(h_{t, \mu}\right)= \zeta_{\mu+\rho_\kk}\left(\xi \right) \cdot d(\mu+\rho_\kk).
\end{align}
where $\zeta_{\mu+\rho_\kk} \colon Z_G\to \C$ denotes the central character of the discrete series representation with Harish-Chandra parameter $\mu + \rho_\kk$, and 
\begin{equation}
\label{formal degree}
d(\mu + \rho_\kk) =   \frac{1}{(2 \pi)^{\frac{\dim \kp}{2}} 2^{(\frac{\dim \kp}{2}-1) / 2}} \cdot \frac{\prod_{\alpha \in R^+\left(\kg, \kt\right)}(\mu + \rho_\kk,  \alpha)}{\prod_{\alpha \in R^+\left(\kk, \kt\right)}\left(\rho_\kk, \alpha\right)}, 	
\end{equation}
is its formal degree.

\begin{proof}
\cite[Section 4.2]{Barbasch83}.
\end{proof}
\end{theorem}

\begin{theorem}[Hyperbolic Orbital Integral]	
\label{hyperbolic}
If $\xi$ is hyperbolic, then 
\begin{align}
\tau_\xi\left(h_{t, \mu}\right) = 0.
\end{align}
\begin{proof}
By \cite{MR0450461}, we know that $\tau_\xi$ has no discrete part contribution for hyperbolic elements. The theorem follows immediately from the Selberg principle in Theorem \ref{cor K selberg}. 
\end{proof}
\end{theorem}

To summarize, for all semisimple elements, we have explicit formulas for the orbital integral in terms of roots. Moreover, Hochs-Wang \cite[Proposition 4.12]{MR3747042} proved an Atiyah-Bott fixed point-type formula:

\begin{theorem}[Topological Formula for Semi-Simple Orbital Integrals]
\label{elliptic II}	
 If  $\xi$ is a semisimple element in $G$, then 
\begin{align}
\tau_\xi\left(h_{t, \mu}\right)= \int_{X^\xi} c^\xi \cdot  \frac{\hat{A}\left(X^\xi\right) \operatorname{Ch}\left(\left[V_\mu\big|_{X^\xi}\right](\xi)\right)}{\operatorname{det}\left(1-\xi e^{-R^{N} / 2 \pi i}\right)^{1 / 2}}.
\end{align}
Here, $c^\xi$ is the cut-off function for the $G_\xi$-action on $X^\xi$, $N$ is the normal bundle of $X^\xi \hookrightarrow X$, and $R^{N}$ is its curvature.	
\end{theorem}

%%%%%%%%%%%%%%%%%%%%%%%%%%%%%%%%%%%%%%%%%%%%%%%%%%%%%%%%%%%%%%%%%%%%%%%%%%%%%%%%%%%%%%%%%%%%%%%%%%%%%%%%%%%%%%%%%%%%%%%%
\subsection{Non-semi-simple orbital integral}

When the lattice $\Gamma$ is not uniform, elements in $\Gamma$ are no longer necessarily semi-simple. Consequently, we must consider orbital integrals of arbitrary elements in the Selberg trace formula, as derived in \cite{MR535763}.

The approach involves expressing the unipotent orbital integral as a limit of semisimple orbital integrals. To make this precise, suppose $u = \exp(X)$ is a unipotent element. Then $X$ can be embedded in a Lie triple $\{X, H, Y\}$ such that $\theta X = -Y$, where $\theta$ is the Cartan involution on $G$. Define
\begin{align}
Z_0 = X - Y \in \kk, \quad \mathfrak{U} = X + Z_Y,
\end{align} 
where $Z_Y$ is the centralizer of $Y$ in $\kg$. It is shown in \cite{MR320232} that $\mathfrak{U}$ is transverse to the adjoint action of $G$ on the Lie algebra $\kg$, and the union $\bigcup_{t>0} O(t(X - Y))$ contains the orbit of $X$ in its closure. This allows the orbital integral for $X$ to be related to the one for $t Z_0$ ($t > 0$), which are orbital integrals of elliptic elements. The Fourier transform of the orbital integral of $u = \exp(X)$ can be found in \cite[Theorem 7.1]{MR526310}. 

For an arbitrary element $\xi \in G$, let 
\begin{align}
\xi = \eta u, \quad u = \exp(X)
\end{align}
be the Jordan decomposition, where $\eta$ is semisimple, $X$ is unipotent, and $s$ commutes with $u$. For $t > 0$, define 
\begin{align}
\xi_t =  z_t, \quad z_t = \operatorname{exp}(t Z_0).
\end{align} 
The orbital integral of $\xi$ can be approximated by the semisimple orbital integral of $\xi_t$. Let $\xi_0$ be the element defined in (\ref{equ Z0}). 
The following theorem is due to Barbasch \cite[Theorem 4.3.2]{MR2626602} or \cite[Theorem 7.1]{MR526310}. 
\begin{theorem}
\label{thm barbasch}
Suppose that $\xi = \eta\exp(X)$. There exists a constant $c_\xi$ such that  
\begin{align}
 \tau_{\xi}(f) =  c_\xi \cdot \sum_{\mu \in \Lambda_\kt^*} \left( \left[\overline{\left\langle \mu, Z_0 \right \rangle}\right]^{\frac{\dim \kn_{\eta,1}}{2}} \cdot \prod_{\beta \in R^+(\xi_0) } \left \langle \mu, \beta\right\rangle \cdot e^{\mu}(\eta) \right) \cdot \Theta_{\mu}(f) + \tau_\xi^{(c)}(f),
\end{align}
where  the $c_\xi$ is the composite of certain other constants, each of which is explicitly computable, except, perhaps, for one, namely, the volume of
\begin{align}
G_\eta \cdot Z_0 \cap\left(X+\kg_{\eta, Y}\right)
\end{align}
in a canonically defined finite measure arising from fiber integration; details can be found in \cite{MR2626602}.

\end{theorem}

\begin{lem}
\label{unipotent su}
The discrete part of the orbital integral $\tau_\xi$ vanishes unless 
\begin{align}
\kg_\eta \cong \mathfrak{su}(n,1) \quad \text{and} \quad X \in \kn_{\eta, 2}.
\end{align}  

\begin{proof}
By \cite[Lemma 8.1]{MR526310}, the orbital integral has no discrete series contribution if there exists a Weyl group element $w \in W_{\kg_\eta}$ such that the Lie triple $\{X, H, Y\}$ is conjugate to $\{-X, H, -Y\}$. 

It is shown in \cite[Theorem 3.2]{MR526310} that this is always the case unless
\begin{align}
\dim \kn_{\eta,2} = 1, \quad \text{and} \quad X \in \kn_{\eta,2}.
\end{align}	
The condition $\dim \kn_{\eta,2} = 1$ implies that $\kg_\eta \cong \mathfrak{su}(n,1)$. This completes the proof.
\end{proof}
\end{lem}
\begin{remark}
\label{remark ceta}
For the case of $\kg_\eta \cong \mathfrak{su}(n,1)$, the constant $c_{\xi = \eta \exp\left(X\right)}$ depends only on the connected component of $X \in \kn_{\eta, 2} \setminus \{0\}$. Thus, we denoted it by $c^\pm_\eta$. One can check that 
\begin{align}
c_\eta^+ = (-1)^{n}\cdot  c_\eta^-. 
\end{align} 	
\end{remark}
Combining Theorem \ref{cor K selberg}, Theorem \ref{thm barbasch} and Lemma \ref{unipotent su}, we obtain the following result:
\begin{cor}[Non-Semi-Simple Orbital Integral]
\label{unipotent}
Suppose that $\xi = \eta\exp(X)$ with $\eta$ semisimple and $X$ unipotent. Then
\begin{align}
\tau_\xi\left(h_{t, \mu}\right)= 0
\end{align}
unless $\kg_\eta \cong \mathfrak{su}(n,1)$ and $X \in \kn_{\eta, 2}$. In this case, let $\{X, H, Y\}$ be the associated Lie triple, and define $Z_0 = X - Y$. For $\xi^\pm = \eta \exp \left(\pm X\right)$,  we have that \begin{align*}
\tau_{\xi^\pm }\left(h_{t, \mu}\right)
=& (-1)^{\frac{\dim \kp}{2}}\cdot c^\pm_\eta \cdot  \sum_{w \in W_\kk} \Bigg( \left[\overline{\left\langle w( \mu+\rho_\kk), Z_0 \right \rangle}\right]^{\frac{\dim \kn_{\eta,1}}{2}} \\
&\times \prod_{\beta \in R^+(\xi_0) } \left\langle w( \mu+\rho_\kk), \beta\right\rangle \cdot e^{w( \mu+\rho_\kk)}(\eta) \Bigg).
\end{align*}
\end{cor}

%%%%%%%%%%%%%%%%%%%%%%%%%%%%%%%%%%%%%%%%%%%%%%%%%%%%%%%%%%%%%%%%%%%%%%%%%%%%%%%%%%%%%%%%%%%%%%%%%%%%%%%%%%%%%%%%%%%%%%%%
\subsection{Weighted orbital integral and weighted character}
\label{sec weighted}

Fix a maximal compact subgroup $K$ of $G$. By the Iwasawa decomposition, one can define maps $H_P, H_{\overline{P}}: G \rightarrow \ka$ as follows:
\begin{align}
H_P(k a n) = \log a, \quad H_{\overline{P}}(k a \overline{n}) = \log a, \quad \forall k \in K, a \in A, n \in N, \overline{n} \in \overline{N}.
\end{align}
Let $\lambda_P$ be the positive multiple of the roots of $\mathfrak{a}$ in $\mathfrak{n}$ satisfying 
\begin{align}
\lambda_P(H) = \pm|H|, \quad H \in \ka. 
\end{align} 
Define the function 
\begin{align}
w(x) = \lambda_P\left(H_P(x) - H_{\overline{P}}(x)\right),
\end{align} 
which is positive-valued and satisfies
\begin{align}
w(x m) = w(x), \quad m \in M.
\end{align}

\begin{definition}
For $\xi \in MA$, the \emph{weighted orbital integral} is defined as
\begin{equation}
\label{weighted orbi}	
\tau^{\text{w}}_\xi(f) = \int_{G / G_\xi} f\left(x \xi x^{-1}\right) \cdot w(x) \; dx, \quad f \in \mathcal{C}(G).
\end{equation}	
The integral is absolutely convergent and defines a non-invariant tempered distribution on $G$. It is independent of $\lambda_P$ and proportional to the Haar measure on $G$. 
\end{definition}

The Fourier transform of the weighted orbital integral has been computed in \cite{MR853551, MR681613, MR412348, MR808270, AM12, MR1461205, MR1282200}.

In Theorem \ref{elliptic I}, we established that the restriction of the character of a discrete series representation to $T^{\text{reg}}$ relates to elliptic orbital integrals. However, its restriction to the non-compact Cartan subgroup $H$ cannot be detected by hyperbolic orbital integrals (Theorem \ref{hyperbolic}). Surprisingly, Arthur proved in \cite{MR0412348} that it is instead related to weighted orbital integrals.

\begin{theorem}
\label{weighted thm}
Suppose that $\xi = ma \in H$. Then,
\begin{equation}
\label{weighted equ}
\tau^{\text{w}}_\xi(f) = -\sum_{\omega \in \widehat{G}_d} \Omega_{\omega}(h) \cdot \Theta_\omega\left(f\right) + \text{continuous part}.
\end{equation}
\begin{proof}
See \cite[Theorem 14.1]{MR808270} or \cite{MR0412348}. 
\end{proof}
\end{theorem}

\begin{remark}
In equation \eqref{weighted equ}, the formula holds for $A$-regular points $\xi = ma$ with $a \neq 1$. Based on considerations from the Selberg trace formula, one may need to take the limit $a \to 1$ and consider $\tau^{\text{w}}_m(f)$. Warner extended equation \eqref{weighted equ} to the case $a = 1$ in \cite[Section 15]{MR808270}.
\end{remark}

Recall the normalized intertwining operator defined in equation (\ref{normal intertwining}):
\begin{align}
U_w \colon H_{\sigma, \nu} \rightarrow H_{\sigma,-\nu}.
\end{align}
If $\phi$ is a meromorphic function of $\nu$, we use the abbreviation
\begin{align}
\frac{\partial}{\partial \nu} \phi(\nu) = \left.\frac{d}{dz} \phi\left(z \cdot \nu\right)\right|_{z=0}.
\end{align}
For $\sigma \in \widehat{M}$ with $\mu(\sigma, 0) \neq 0$ and $f \in \mathscr{C}(G)$, we define:

\begin{definition}
The \emph{weighted character} is given by
\begin{align}
\widetilde{\Theta}_{\sigma, \nu}(f) = -\operatorname{Tr} \left(\pi_{\sigma, \nu}(f) \circ U(\sigma, \nu)^{-1} \circ \frac{\partial}{\partial \nu} U(\sigma, \nu) \right),
\end{align}
where the derivative is taken for the operator in $H_{\sigma, \nu}$. The weighted character $\widetilde{\Theta}$ is a non-invariant tempered distribution, proportional to the Haar measure on $G$ and to $\lambda_P$. It does not depend on the Haar measure on $\overline{N}$ used to define $U(\sigma, \nu)$.	
\end{definition}

In \cite[Section 10]{MR625344}, Arthur defined invariant distributions. In our real rank one setting, this procedure simplifies significantly.

\begin{definition}
For any $\xi \in L$ and $f \in \mathcal{C}(G)$, we define an invariant tempered distribution by the formula:
\begin{align}
\widetilde{W}(\xi, f) = \tau^{\text{w}}_\xi(f) - \frac{1}{2 \pi i} \sum_{\sigma \in \widehat{M}} \int_{i\ka^*} \overline{\Theta_{\sigma, \nu}(\xi)} \cdot\widetilde{\Theta}_{\sigma, \nu}(f) \; d\nu.
\end{align}
In particular, this induces a map:
\begin{align}
\widetilde{W}_\xi \colon K_*\left(C^*_r G\right) \to \C. 
\end{align}
\end{definition}

\begin{cor}
\label{invariant thm}
For any $\xi \in L$,
\begin{align}
\widetilde{W}_\xi \left(h_{t, \mu}\right) =  (-1)^{\frac{\dim \kp}{2}+1} \Omega_{\mu + \rho_\kk}(\xi).
\end{align}	
\begin{proof}
By conjugation, we may assume that $\xi \in H$, where $H$ is the non-compact Cartan subgroup. By definition, the invariant distribution $\widetilde{W}_\xi$ and the non-invariant distribution $\tau^{\text{w}}_\xi$ have the same discrete part. The corollary follows from Theorem \ref{weighted thm} and  \ref{cor K selberg}.
\end{proof}
\end{cor}

\begin{remark}
\label{IL rem}
It has been shown in \cite[Lemma 5]{MR671316} that
\begin{align}
\Omega_{\mu + \rho_\kk}(e) = 0, \quad \text{unless } \kg = \mathfrak{sl}(2, \R) \text{ or } \mathfrak{su}(2, 1).  
\end{align}	
\end{remark}
%%%%%%%%%%%%%%%%%%%%%%%%%%%%%%%%%%%%%%%%%%%%%%%%%%%%%%%%%%%%%%%%%%%%%%%%%%%%%%%%%%%%%%%%%%%%%%%%%%%%%%%%%%%%%%%%%%%%%%%%
%%%%%%%%%%%%%%%%%%%%%%%%%%%%%%%%%%%%%%%%%%%%%%%%%%%%%%%%%%%%%%%%%%%%%%%%%%%%%%%%%%%%%%%%%%%%%%%%%%%%%%%%%%%%%%%%%%%%%%%%

\section{Trace formula and proof of main results}
\label{sec trace}
In this section, we prove the main result. The primary tool we use is the invariant trace formula for Hecke operators. The general case for the adelic setting was developed by Arthur in \cite{MR625344, MR928262, MR939691}. In this paper, we employ the trace formula for real groups, which was established by Hoffmann \cite{MR1001467, MR1724290}. It is worth noting that the specific formula we use is simpler than Hoffmann's general result, as we impose the additional assumption that $G$ has real rank one.

We assume that $G$ is a semi-simple Lie group of real rank one. Take $\alpha \in C_G(\Gamma)$ and let $\Xi = \Gamma \alpha \Gamma$ be a double coset. 

 \subsection{Trace Formula for Hecke Operators}

Recall the decomposition of $\Xi$:
\begin{align}
\Xi = \underbrace{\Xi_S \sqcup \Xi_P(r) \sqcup \Xi_P(s, *)}_{\text{semi-simple part}} \sqcup  \Xi_P(s, **).
\end{align}
We begin with the following lemma. 

\begin{lem}
\label{hyper lem}
If $G$ is of real rank one and equal rank, then 
\begin{align}
 \Xi_P(r) \subseteq \Xi_H.
\end{align}
That is, every element in $\Xi_P(r)$ is hyperbolic. We split
\begin{align}
\Xi_H =  \Xi_P(r) \sqcup \Xi^s_H, \quad \Xi^s_H  = \Xi_P \setminus  \Xi_P(r). 
\end{align}
\begin{proof}
Take $\xi \in \Xi_P(r)$. Since $\xi$ is semi-simple, it follows from \cite[Lemma 7.10]{MR1724290} that $\xi$ is $N$-conjugate to an element in $L$. Let $\eta \in \{\xi\}_N \cap L \neq \emptyset$. We know that $\eta$ is regular in the sense that $N_\eta = \{1\}$. The assumption that $G$ has real rank one and equal rank excludes the case $\kg \cong \mathfrak{so}(2n+1, 1)$. It has been shown in \cite[Page 26]{MR614517} that the centralizers $N_m$ are non-trivial for all $m \in M$ when $\kg \neq \mathfrak{so}(2n+1, 1)$. Thus,  
\begin{align}
\eta = ma \in L = MA,
\end{align}
with $a \neq 1$ in $A$. Consequently, $\xi$ can be conjugated to an element in the non-compact Cartan subgroup and must be hyperbolic.
\end{proof}
\end{lem}

\begin{theorem}[Non-Invariant Trace Formula]
\label{non invariant thm}
For any $f \in\mathcal{C}^1(G)$, the trace of $T_\Xi \circ R^{\Gamma}_d(f)$ equals the sum of the following parts:
\begin{enumerate}
    \item \textbf{Semi-simple orbital integrals in} $\Xi_C \sqcup \Xi_E \sqcup \Xi^s_H$: 	
    \begin{align*}
    C_{\Xi}(f) &= \operatorname{vol}\left(\Gamma \bs G\right) \sum_{\{\xi\}_{\Gamma} \subset \Xi_C} \tau_\xi(f),\\
    E_{\Xi}(f) &= \sum_{\{\xi\}_{\Gamma} \subset \Xi_E} \operatorname{vol}\left(\Gamma_{\xi} \bs G_{\xi} \right) \cdot  \tau_\xi(f),\\
    H_{\Xi}(f) &= \sum_{\{\xi\}_{\Gamma} \subset \Xi^s_H} \operatorname{vol}\left( \Gamma_{\xi} \bs G_{\xi} \right) \cdot  \tau_\xi(f).
    \end{align*}
    
    \item \textbf{Semi-simple orbital integrals from} $\Xi_P(r)$:
    \begin{align}
    R_\Xi(f) = \frac{1}{2}  \sum_{P \in E(G, \Gamma)}  \sum_{\{\xi\}_{\Gamma_M} \subset \Xi_L(\mathrm{r})} \alpha^\Xi_\xi \cdot \tau_\xi(f),
    \end{align}
    where $\alpha^\Xi_\xi$ are some constants defined in \cite[Page 105]{MR1724290}. 
    
    \item \textbf{Unipotent orbital integrals}:
    \begin{align}
    U_{\Xi}(f) = \sum_{P \in E(G, \Gamma)} \sum_{\{\eta\}_{\Gamma_M} \subset \Xi_{L}(s)}  \sum_{\left\{X\right\}_{L_\eta} \subset \overset{\circ}{N}_{\eta, 1} \cup \overset{\circ}{N}_{\eta, 2}}C_P(\Xi, s, X) \cdot\tau_{\eta\exp\left (X\right)}(f),
    \end{align}
    where $C_P(\Xi, \eta, X)$ are defined as the constant or residue terms of certain Zeta functions associated with $\Xi$.
    
    \item \textbf{Weighted orbital integrals}:
    \begin{align*}
    W_{\Xi}(f) = &\frac{1}{2}  \sum_{P \in E(G, \Gamma)} \sum_{\{\eta\}_{\Gamma_M} \subset\Xi_{L}} \operatorname{vol}\left(M_\eta /\left(\Gamma_M\right)_\eta\right) \\
    &\times \left|\operatorname{det} \operatorname{Ad}_{\kn}(\eta)\right|^{-1 / 2}  \cdot\#\left[ \Xi \cap \eta N \colon \Gamma \cap N \right] \cdot  \tau^{\text{w}}_\eta(f).
    \end{align*}
    
    \item \textbf{Intertwining term}:
    \begin{align}
    I_{\Xi}(f) = \frac{1}{4 \pi i} \sum_{\substack{\sigma \in \widehat{M} \\ \sigma \mid_{Z_{\Gamma}}= \mathrm{id}}}\int_{i\ka^*} \Tr\left(\pi_{\sigma, v}(f) \circ T_{\alpha, \pi_{\sigma, \nu}}  \circ U^\Gamma\left(\sigma,-\nu\right) \circ \frac{d}{d \nu} U^\Gamma\left(\sigma,  \nu\right)\right) \; d \nu.
    \end{align}
    
    \item \textbf{Residue term}:
    \begin{align}
    J_{\Xi}(f)  = -\frac{1}{4}  \sum_{\substack{\sigma \in \widehat{M} \\ \sigma \mid_{Z_{\Gamma}}= \mathrm{id}}}\Tr \left(\pi_{\sigma, 0}(f) \circ  T_{\alpha, \pi_{\sigma, 0}}  \circ U^\Gamma\left(\sigma,  0\right)\right).
    \end{align}
\end{enumerate}
\end{theorem}

\begin{remark}
In \cite{MR1001467, MR1724290}, the trace formula in Theorem \ref{non invariant thm} does not require $\Gamma$ to be arithmetic but was originally stated only for functions $f \in C_c(G)$. Under the assumption that $\Gamma$ is arithmetic, the non-adelic trace formula can be seen as a special case of the adelic trace formula. In \cite{MR2811597, MR3534542, MR2434856}, the adelic trace formula has been extended to Harish-Chandra $L^1$-Schwartz classes. In fact, if Hoffmann's formula could be extended to a broader class of functions with non-compact support (which seems to be feasible by similar argument in \cite[Section 8]{MR535763}), it might be possible to drop the assumption that $\Gamma$ is arithmetic.

\end{remark}

In \cite{Barbasch83}, Barbasch and Moscovici applied the special case of Theorem \ref{non invariant thm} to the function 
\begin{align}
h_{t, \mu}=\Tr^{E_\mu^+}\left( h_{t, \mu}^{+} \right)-\Tr^{E_\mu^-}\left( h_{t, \mu}^{-}\right).
\end{align}
That is, they obtained the formula:
\begin{align*}
\Tr\left(R_d^{\Gamma}\left(h_{t, \mu}\right) \right) =& C_{\Gamma}(h_{t, \mu})  + E_{\Gamma}(h_{t, \mu})  + H_{\Gamma}(h_{t, \mu})  + R_{\Gamma}(h_{t, \mu}) 	\\
&+ U_{\Gamma}(h_{t, \mu}) + W_{\Gamma}(h_{t, \mu}) + I_{\Gamma}(h_{t, \mu})  + J_{\Gamma}(h_{t, \mu}).
\end{align*}

Among the terms above, the most challenging to analyze are the \emph{weighted orbital integral} term $W_{\Gamma}(h_{t, \mu})$ and the \emph{intertwining term} $I_{\Gamma}(h_{t, \mu})$, as they are not invariant. In particular, the Selberg principle (Theorem \ref{cor K selberg}) does not apply to these two terms. To address this, Barbasch and Moscovici decomposed $h_{t, \mu}$ as follows:
\begin{align}
h_{t, \mu}=h_{t, \mu}^{(d)}+h_{t, \mu}^{(c)},
\end{align}
where $h_{t, \mu}^{(d)}$ and $h_{t, \mu}^{(c)}$ represent the \emph{discrete} and \emph{continuous} parts of the Plancherel expansion for $h_{t, \mu}$, respectively. Explicitly, these components are given by:
\begin{align}
h_t^{(d)}(x)=\sum_{\omega \in \widehat{G}_d} d(\omega) \cdot  \operatorname{Tr}  \left(\omega\left(h_t\right) \omega(x)^{-1}\right), \quad x \in G,
\end{align}
where $d(\omega)$ denotes the formal degree of $\omega$; and
\begin{align}
h_t^{(c)}(x)=\sum_{\sigma \in \widehat{M}} \int_{i\R} \operatorname{Tr}\left(\pi_{\sigma, v}\left(h_t\right) \pi_{\sigma, v}(x)^{-1}\right) \mu_\sigma(v) d v, \quad x \in G,
\end{align}
where $\mu_\sigma$ is the Plancherel density.

One of the main observations in \cite{Barbasch83} is that if $\mu$ is regular, then the limits
\begin{equation}
\label{vanishing equ}	
\lim_{t\to \infty}  W_{\Xi}\left(h^{(c)}_{t, \mu}\right)  = 0, \quad \lim_{t\to \infty}  I_{\Xi}\left(h^{(c)}_{t, \mu}\right)   = 0
\end{equation}
hold, which leads to the final formula for $\Tr\left(R_d^{\Gamma}\left(h_{t, \mu}\right) \right)$.

However, we emphasize that if $\mu$ is singular or if $\Xi \neq \Gamma$, the limits
\begin{align}
 W_{\Xi}\left(h^{(c)}_{t, \mu}\right), \quad \text{and} \quad I_{\Xi}\left(h^{(c)}_{t, \mu}\right)  
 \end{align}
might not be zero, and their explicit computation remains challenging. To overcome this difficulty, we will employ the following formula. 
\begin{theorem}[Invariant Trace Formula for Hecke Operators]
\label{invariant thm}
Under the same setting as Theorem \ref{non invariant thm}, the following formula holds:
\begin{align*}
\Tr\left(R_d^{\Gamma}\left(f\right) \right) =& C_{\Xi}(f)  + E_{\Xi}(f)  + H_{\Xi}(f)  + R_{\Xi}(f) 	\\
+& U_{\Xi}(f) + \widetilde{W}_{\Xi}(f) + \widetilde{I}_{\Xi}(f)  + J_{\Xi}(f).
\end{align*}
where
\begin{align*}
\widetilde{W}_{\Xi}(f) &= \frac{1}{2} \sum_{P \in E(G, \Gamma)} \sum_{\{\eta\}_{\Gamma_M} \subset\Xi_{L}} \operatorname{vol}\left(M_\eta /\left(\Gamma_M\right)_\eta\right)\\
&\times \left|\operatorname{det} \operatorname{Ad}_{\kn}(\eta)\right|^{-1 / 2} \cdot\#\left[ \Xi \cap \eta N \colon \Gamma \cap N \right] \cdot  \widetilde{W}_\eta(f),
\end{align*}
and
\begin{align}
\widetilde{I}_{\Xi}(f) = \frac{1}{4 \pi i} \sum_{\substack{\sigma \in \widehat{M} \\ \sigma \mid_{Z_{\Gamma}}= \mathrm{id}}} \int_{i \ka^*}\Tr\left( T_{\alpha, \pi_{\sigma, \nu}}  \circ \mathbf{S}\left(\sigma, -\nu \right) \circ \frac{d}{d \nu } \mathbf{S}\left(\sigma, \nu\right)\right)\cdot  \Theta_{\sigma, \nu}(f)\; d \nu.
\end{align}
Here  $ \mathbf{S}\left(\sigma, \nu\right)$ are some  operators defined in \cite[Section 6]{MR1724290}.
\end{theorem}
 
\begin{remark}
The key advantage of Theorem \ref{invariant thm} is that each term in the trace formula is now invariant, making the analysis more structured. However, the computation of the operator $\mathbf{S}(\sigma, \nu)$ remains intricate and less explicit. Despite this complexity, we note that the contribution of $\widetilde{I}_{\Xi}$ vanishes due to Theorem \ref{cor K selberg}, which simplifies our computations significantly.
\end{remark}

%%%%%%%%%%%%%%%%%%%%%%%%%%%%%%%%%%%%%%%%%%%%%%%%%%%%%%%%%%%%%%%%%%%%%%%%%%%%%%%%%%%%%%%%%%%%%%%%%%%%%%%%%%%%%%%%%%%%%%%%

\subsection{Proof of Theorem \ref{main thm I} and \ref{main thm II}}
Applying Theorem \ref{invariant thm} to the function $h_{t, \mu}$, 
\begin{align*}
 L\left(T_\Xi, D^{\Gamma}_\mu\right) =&\Tr\left( T_\Xi\circ  R_d^{\Gamma}\left(h_{t, \mu}\right) \right)\\
 =& C_{\Xi}(h_{t, \mu})  + E_{\Xi}(h_{t, \mu})  + H_{\Xi}(h_{t, \mu})  + R_{\Xi}(h_{t, \mu}) 	\\
+& U_{\Xi}(h_{t, \mu}) + \widetilde{W}_{\Xi}(h_{t, \mu}) + \widetilde{I}_{\Xi}(h_{t, \mu})  + J_{\Xi}(h_{t, \mu}).
\end{align*}
We  compute the above terms separately:
\begin{enumerate}
	\item By Theorem \ref{central}, 
	\begin{align}
	C_{\Xi}(h_{t, \mu})  = 		 \delta(\mu, \Xi) \cdot \operatorname{vol}(\Gamma \backslash G) \cdot d(\mu+\rho_\kk),
	\end{align}
	which vanishes if $\mu$ is singular. 
	\item  By Theorem \ref{elliptic I}, 
	\begin{align*}
	E_{\Xi}(h_{t, \mu}) = &(-1)^{\frac{\dim \kp}{2}}  \sum_{ \left\{\xi\right\}_\Gamma \in \Xi_E} \operatorname{vol}\left(\Gamma_\xi \backslash G_\xi\right) \cdot d_\xi^{-1} \\
        &\times \sum_{w \in W_\kk/W_{\kk_\xi}} \operatorname{det}(w) \cdot\left(\frac{ \prod_{\alpha \in R^+(\xi)}\left\langle w\left(\mu+ \rho_\kk\right), \alpha\right\rangle \cdot  e^{w\left(\mu+\rho_\kk\right)}(\xi)}{ e^{\rho_\kg}\left(\xi\right)\prod_{\beta \in R^+ \setminus R^+(\xi)} \left(1  - e^{-\beta}\left(\xi\right)\right)}\right). 
	\end{align*}
	\item By Theorem \ref{hyperbolic}
	\begin{align}
	 H_\Xi(h_{t, \mu}) = 0. 
	\end{align}
	\item By Theorem \ref{hyperbolic} and Lemma \ref{hyper lem}, 
	\begin{align}
	 R_\Xi(h_{t, \mu}) = 0. 
	\end{align}
	\item By Corollary \ref{unipotent}, 
	\begin{align*}
	U_\Xi(h_{t, \mu})=& \sum_{P \in E(G, \Gamma)} \sum_{\{s\}_{\Gamma_M} \subset \Xi_{L}(s)} (-1)^{r + r_s + \frac{\dim \kp}{2}}\\ &\times \delta(s) \cdot \left\{ C_P(\Xi, s, X) \cdot c^{-1}_{s, X} + C_P(\Xi, s, -X) \cdot c^{-1}_{s, -X} \right\} \\
        &\times \sum_{w \in W_\kk} \left( \left[\overline{\left\langle w( \mu+\rho_\kk), Z_0 \right \rangle}\right]^{\frac{\dim \kn_{s,1}}{2}}\cdot \prod_{\beta \in R^+(\xi_0) }\left \langle w( \mu+\rho_\kk), \beta\right\rangle \cdot e^{w( \mu+\rho_\kk)}\left(s\right) \right)	
	\end{align*}
	\item By Corollary \ref{invariant thm}, 
	\begin{align*}
	\widetilde{W}_\Xi(h_{t, \mu}) = &  \frac{(-1)^{\frac{\dim \kp}{2}}}{2} \sum_{P \in E(G, \Gamma)} \sum_{\{s\}_{\Gamma_M} \subset\Xi_{L}} \operatorname{vol}\left(M_s /\left(\Gamma_M\right)_s\right) \cdot \left|\operatorname{det} \operatorname{Ad}_{\kn}(s)\right|^{-1 / 2}	\\
	\times & \#\left[ \Xi \cap sN \colon \Gamma \cap N \right] \cdot  \Omega_{\mu + \rho_\kk}(s), 
	\end{align*}
	which also vanishes if $\mu$ is singular. 
	\item By Theorem \ref{cor K selberg}, 
	\begin{align}
	\widetilde{I}_{\Xi}(h_{t, \mu}) = 0.  
	\end{align}
	\item If $\mu$ is regular, then it follows from Theorem \ref{cor K selberg} that 
	\begin{align}
	J_{\Xi}(h_{t, \mu}) = 0. 
	\end{align}
	Otherwise, 
	\begin{align}
	J_{\Xi}(h_{t, \mu}) = -\frac{1}{2} \operatorname{Tr}\left(  T_{\alpha, \pi_{\sigma, 0}} \circ U^{\Gamma, +}\left(\sigma,  0\right) \right).
	 \end{align}
\end{enumerate}

\bibliographystyle{alpha}
\bibliography{mybib}

\end{document}